\documentclass{mcom-l}

\numberwithin{equation}{section}

\usepackage{amsmath,amssymb,amsthm}
\usepackage{enumerate}
\usepackage{cases}
\usepackage{color}
\usepackage[mathscr]{euscript}
\usepackage{url}
\usepackage[all]{xy}
\usepackage{booktabs}
\usepackage{graphicx}
\usepackage{epstopdf}
\usepackage{algorithmic}
\usepackage[caption=false]{subfig}
\usepackage{mathtools}
\usepackage{soul} 
\usepackage{multirow}
\usepackage{makecell}
\usepackage{bigdelim}
\usepackage[ruled,vlined]{algorithm2e}
\usepackage{threeparttable}
\usepackage{hyperref}
\usepackage{xcolor}
\definecolor{cobalt}{rgb}{0.0, 0.28, 0.67}
\definecolor{darkblue}{rgb}{0.0, 0.0, 0.55}

\hypersetup
{
bookmarksopen=true,
pdftitle={Nonsymmetric Jacobi-G algorithm},
pdfauthor={SL},
pdfsubject={15A12},
pdfmenubar=true,
pdfhighlight=/O,
colorlinks=true,
pdfpagelayout=SinglePage,
pdffitwindow=true,
linkcolor=cobalt,
citecolor=cobalt,
urlcolor=cobalt
}

\usepackage[capitalize,nameinlink,noabbrev]{cleveref}

\usepackage[scale=0.75]{geometry}

\theoremstyle{definition}

\newtheorem{theorem}{Theorem}[section]
\newtheorem{corollary}[theorem]{Corollary}
\newtheorem{lemma}[theorem]{Lemma}
\newtheorem{definition}[theorem]{Definition}
\newtheorem{proposition}[theorem]{Proposition}
\newtheorem{remark}[theorem]{Remark}

\newtheorem{example}[theorem]{Example}

\definecolor{darkred}{rgb}{0.7,0,0}
\definecolor{darkgreen}{rgb}{0,0.46,0}
\definecolor{purple}{rgb}{0.6,0,0.5}
\newcommand{\fin}{\color{black}}

\newcommand{\tens}[1]{\boldsymbol{\mathcal{#1}}}
\newcommand{\tenselem}[1]{\mathcal{#1}}

\newcommand{\matr}[1]{\boldsymbol{#1}}

\newcommand{\vect}[1]{\boldsymbol{#1}}

\newcommand{\set}[1]{\mathscr{#1}}
\newcommand{\T}{{\sf T}}    
\renewcommand{\H}{{\sf H}}   

\makeatletter
\newcommand{\rank}[1]{\mathop{\operator@font rank}(#1)}
\newcommand{\colrank}[1]{\mathop{\operator@font colrank}\{#1\}}
\newcommand{\krank}[1]{\mathop{\operator@font krank}\{#1\}}
\newcommand{\trace}[1]{\mathop{\operator@font tr}\left(#1\right)}
\newcommand{\symmm}[1]{\mathop{\operator@font sym}\left(#1\right)}
\newcommand{\skeww}[1]{\mathop{\operator@font skew}(#1)}
\newcommand{\Diag}[1]{\mathop{\operator@font Diag}\{#1\}}  
\newcommand{\diag}[1]{\mathop{\operator@font diag}\{#1\}}  
\newcommand{\Span}[1]{\mathop{\operator@font Span}\{#1\}}  
\newcommand{\argmin}{\mathop{\operator@font argmin}}

\newcommand{\offdiag}[1]{\mathop{\operator@font offdiag}\{#1\}}  
\newcommand{\Proj}[2]{\mathop{\operator@font Proj_{#1}}{#2}}
\newcommand{\ProjGrad}[2]{\mathop{{\operator@font grad} }#1(#2)}

\newcommand{\expp}[1]{\mathop{\operator@font exp}\left(#1\right)}
\makeatother

\newcommand{\eqdef}{\stackrel{\sf def}{=}}
\newcommand{\RR}{\mathbb{R}}
\newcommand{\CC}{\mathbb{C}}

\newcommand{\contr}[1]{\mathop{\bullet_{#1}}}  

\newcommand{\Utwo}{\Psi}

\newcommand{\RInner}[2]{\left\langle #1, #2\right\rangle_{\Re}}  
\newcommand{\R}{\Re}
\newcommand{\I}{\Im}
\newcommand{\ui}{\mathbf{i}}

\newcommand{\p}[1]{\mathop{#1}\limits_{(p)}}
\newcommand{\UN}[1]{\set{U}_{#1}}
\newcommand{\Gmat}[3]{\matr{G}(#1,#2,#3)}

\newcommand{\ue}{\mathrm{e}}

\newcommand{\hij}[3]{h_{(#1,#2), #3}}

\newcommand{\hijtilde}[3]{\tilde{h}_{(#1,#2), #3}}

\makeatletter
\def\bbordermatrix#1{\begingroup \m@th
\@tempdima 4.75\p@
\setbox\z@\vbox{%
\def\cr{\crcr\noalign{\kern2\p@\global\let\cr\endline}}%
\ialign{$##$\hfil\kern2\p@\kern\@tempdima&\thinspace\hfil$##$\hfil
&&\quad\hfil$##$\hfil\crcr
\omit\strut\hfil\crcr\noalign{\kern-\baselineskip}%
#1\crcr\omit\strut\cr}}%
\setbox\tw@\vbox{\unvcopy\z@\global\setbox\@ne\lastbox}%
\setbox\tw@\hbox{\unhbox\@ne\unskip\global\setbox\@ne\lastbox}%
\setbox\tw@\hbox{$\kern\wd\@ne\kern-\@tempdima\left[\kern-\wd\@ne
\global\setbox\@ne\vbox{\box\@ne\kern2\p@}%
\vcenter{\kern-\ht\@ne\unvbox\z@\kern-\baselineskip}\,\right]$}%
\null\;\vbox{\kern\ht\@ne\box\tw@}\endgroup}
\makeatother

\begin{document}

\title[Jacobi-type algorithms for homogeneous polynomial optimization]{Jacobi-type algorithms for homogeneous polynomial optimization on Stiefel manifolds with applications to tensor approximations}

\author{Zhou Sheng}
\address{Department of Data Science, Anhui University of Technology, Maanshan, Anhui, China}
\email{szhou03@live.com}
\thanks{The first author was supported in part by the Anhui Provincial Natural Science Foundation (No. 2208085QA07) and the Youth Foundation of Anhui University of Technology (No. QZ202114).}

\author{Jianze Li}
\address{Shenzhen Research Institute of Big Data, The Chinese University of Hong Kong, Shenzhen, Guangdong, China}
\email{lijianze@gmail.com}
\thanks{The second author was supported in part by the National Natural Science Foundation of China (No. 11601371) and the GuangDong Basic and Applied Basic Research Foundation (No. 2021A1515010232).}

\author{Qin Ni}
\address{Department of Mathematics, Nanjing University of Aeronautics and Astronautics, Nanjing, Jiangsu, China}
\email{niqfs@nuaa.edu.cn}
\thanks{The third author was supported by the National Natural Science Foundation of China (No. 11771210).}

\subjclass[2020]{Primary 15A69, 90C23; Secondary 65F99, 90C30}

\date{}


\keywords{Homogeneous polynomial optimization, Stiefel manifold, joint approximate tensor diagonalization, 
Jacobi-type algorithm, convergence analysis, \L{}ojasiewicz gradient inequality}

\begin{abstract}
 This paper mainly studies the gradient-based Jacobi-type algorithms to maximize two classes of homogeneous polynomials with orthogonality constraints, and establish their convergence properties. For the first class of homogeneous polynomials subject to a constraint on a Stiefel manifold, we reformulate it as an optimization problem on a unitary group, which makes it possible to apply the gradient-based Jacobi-type (Jacobi-G) algorithm. Then, if the subproblem can always be represented as a quadratic form, we establish the global convergence of Jacobi-G under any one of  three conditions. The convergence result for the first condition is an easy extension of the result in [Usevich et al. SIOPT 2020], while other two conditions are new ones. This algorithm and the convergence properties apply to the well-known joint approximate symmetric tensor diagonalization. 
For the second class of homogeneous polynomials subject to constraints on the product of Stiefel manifolds, we reformulate it as an optimization problem on the product of unitary groups, and then develop a new gradient-based multi-block Jacobi-type (Jacobi-MG) algorithm to solve it. We establish the global convergence of Jacobi-MG under any one of the above three conditions, if the subproblem can always be represented as a quadratic form. This algorithm and the convergence properties are suitable to the well-known joint approximate tensor diagonalization. 
As the proximal variants of Jacobi-G and Jacobi-MG, we also propose the Jacobi-GP and Jacobi-MGP algorithms, and establish their global convergence without any further condition. 
Some numerical results are provided indicating the efficiency of the proposed algorithms.
\end{abstract}

\maketitle

\section{Introduction}\label{sec:introduction}

\subsection{Notation}
Let $\tens{A}\in\CC^{n_1\times n_2\times\cdots\times n_d}$ be an order-$d$ complex tensor, and $\matr{X}\in\CC^{r\times n_p}$ be a complex matrix. We use the \emph{$p$-mode product} defined as 
$(\tens{A}\contr{p}\matr{X})_{q_1\cdots j\cdots q_d}\eqdef\sum_{q_p=1}^{n_p} \tenselem{A}_{q_1\cdots q_p\cdots q_d} X_{j q_p}.$ 
For a matrix $\matr{X}\in\CC^{n\times r}$, $\matr{X}^{\T}$, $\matr{X}^{*}$ and $\matr{X}^{\H}$ represent its \emph{transpose, conjugate} and \emph{conjugate transpose,} respectively.
 In this paper, we will write $(\cdot)^{\dagger}$ and it can stand for $(\cdot)^{\H}$ or $(\cdot)^{\T}$, depending on the context.
The notation $\diag{\tens{A}}$ refers to the vector formed by the diagonal elements of a tensor $\tens{A}\in\CC^{n_1\times n_2\times\cdots\times n_d}$, and it is defined as
$\diag{\tens{A}}\eqdef (\tenselem{A}_{11\cdots 1}, \tenselem{A}_{22\cdots 2}, \cdots, \tenselem{A}_{n_{\min}n_{\min}\cdots n_{\min}})^{\T}$ with $n_{\min} = \min_{1\leq p\leq d}n_p$. We say that a tensor $\tens{A}$ is \emph{diagonal}, if its elements except the diagonal elements are all zero, that is, $\tenselem{A}_{q_1q_2\cdots q_d}=0$ unless $q_1=q_2=\cdots =q_d$ for $1\leq q_p \leq n_p$ and $1\leq p\leq d$. 
We denote by $\|\cdot\|$ the Frobenius norm of a tensor or matrix, and the Euclidean norm of a vector. 
We denote by $\tens{A}^{(r_1,r_2,\cdots,r_d)}$ the $(r_1,r_2,\cdots,r_d)$-dimensional subtensor obtained from  $\tens{A}$ by allowing its indices $q_p$ to take values only in $1\leq q_p\leq r_p$ for $1\leq p\leq d$.

Let 
$\text{St}(r,n,\CC) \eqdef\{\matr{X}\in\CC^{n\times r}: \matr{X}^{\H}\matr{X}=\matr{I}_r\}$
be the \emph{Stiefel manifold} with $1\leq r\leq n$. 
We denote by  
$\Omega\eqdef\text{St}(r_1,n_1,\CC)\times\text{St}(r_2,n_2,\CC)\times\cdots\times\text{St}(r_d,n_d,\CC)$
the product of $d$ Stiefel manifolds $(d>1)$, and
$\Omega^{(p)}\eqdef\text{St}(r_1,n_1,\CC)\times\cdots\times\text{St}(r_{p-1},n_{p-1},\CC)\times\text{St}(r_{p+1},n_{p+1},\CC)\times\cdots\times\text{St}(r_d,n_d,\CC).$ 
Denote that $\vect{\omega}=(\matr{X}^{(1)}, \matr{X}^{(2)},\cdots,\matr{X}^{(d)})\in\Omega$ and $\matr{X}=(\vect{x}_1,\vect{x}_2,\cdots,\vect{x}_r)\in\CC^{n\times r}$.
 Let $\tens{B}\in\CC^{n\times n\times\cdots\times n}$ be an order-$2\kappa$ complex tensor with $\kappa\geq 1$ and $\matr{X}\in\text{St}(r,n,\CC)$.
We denote that 
\begin{equation*}
\matr{\Gamma}(\tens{B},\matr{X},\kappa)\eqdef\tens{B} \contr{1} \matr{X}^\H\contr{2}\cdots \contr{\kappa} \matr{X}^\H\contr{\kappa+1} \matr{X}^\T\contr{\kappa+2}\cdots \contr{2\kappa} \matr{X}^\T. 
\end{equation*}

\subsection{Problem formulation}
This paper mainly considers two important classes of homogeneous polynomial optimization problems with orthogonality constraints. The first class can be written as the \emph{homogeneous polynomial optimization on a Stiefel manifold} (HPOSM): 
\begin{align}
\max\ g:\ \text{St}(r,n,\CC) \longrightarrow \RR,\ 
\matr{X}\longmapsto \sum_{q=1}^{r}
\matr{\Gamma}(\tens{B}^{(q)},\vect{x}_{q},\kappa), 
\label{definition-f-new}
\end{align}
 where $\tens{B}^{(q)}\in\CC^{n\times n\times\cdots\times n}$ is an order-$2\kappa$ complex tensor with $1\leq q\leq r$ and $\kappa\geq 1$.
For $1\leq q\leq r$, we always assume that the tensor $\tens{B}^{(q)}$ is \emph{Hermitian} \cite{nie2020hermitian,ULC2019}, \emph{i.e.},
$\tenselem{B}^{(q)}_{i_1,\cdots,i_\kappa,i_{\kappa+1},\cdots,i_{2\kappa}} = (\tenselem{B}^{(q)}_{i_{\kappa+1},\cdots,i_{2\kappa},i_1,\cdots,i_\kappa})^{*}.$ 
 The second class can be written as the \emph{homogeneous polynomial optimization on the product of Stiefel manifolds} (HPOSM-P): 
\begin{equation}\label{definition-f-prob-1}
\max_{\vect{\omega}\in\Omega}f(\vect{\omega}),
\end{equation}
 where, for each $1\leq p\leq d$, the restricted function
\begin{align}
f_{(p)}:\ \text{St}(r_p,n_p,\CC) \longrightarrow \RR, \ 
\matr{X} \longmapsto f(\vect{\nu}^{(p)},\matr{X})
\eqdef f(\matr{X}^{(1)}, \cdots, \matr{X}^{(p-1)}, \matr{X}, \matr{X}^{(p+1)}, \cdots, \matr{X}^{(d)}),\label{definition-h}
\end{align}
is of form \eqref{definition-f-new} for any fixed $\vect{\nu}^{(p)}\in\Omega^{(p)}$.

\subsection{Homogeneous polynomial optimization with orthogonality constraints}
The homogeneous polynomial optimization has been widely used in applied and computational mathematics. 
In general, there exist two popular tools to handle it. 
The first one is to transform it into a convex optimization problem which can be solved in polynomial time by the \emph{semi-definite programming} (SDP) relaxation. The other one is the \emph{moment sums of squares} (moment-SOS) relaxation. 
For instance, Luo and Zhang \cite{Luo2010A} presented a general semi-definite relaxation scheme for general quartic polynomial optimization under homogeneous quadratic constraints, which leads to a quadratic optimization problem with linear constraints over the semi-definite matrix cone. Ling, Nie, Qi and Ye \cite{Ling2009Biquadratic} presented a bi-linear SDP relaxation for bi-quadratic polynomial optimization over unit spheres. 
Except for these relaxation techniques, He, Li and Zhang \cite{He2010Approximation} also proposed an approximation algorithm for the homogeneous polynomial optimization with quadratic constraints.  
One can refer to \cite{Li2012Approximation} for a survey on the approximation methods.

In the last three decades, the homogeneous polynomial optimization with orthogonality constraints has also received a lot of attention because of the wide applications \cite{edelman1998geometry,Jiang2015A,Wen2013A}, including the \emph{linear eigenvalue problem} \cite{GoluV96:jhu}, \emph{independent component analysis} \cite{Cardoso93:JADE,Como94:sp,Como10:book} and  \emph{orthogonal Procrustes problem} \cite{elden1999procrustes}. 
In general, it is difficult to solve due to the nonconvexity of orthogonality constraints, which may result in not being guaranteed to find a global solution. 
From the perspective of Euclidean optimization, it is a constrained optimization problem, and thus there are many algorithms to solve it; see \cite{Gao2018A,Lai2014A,Hu2020An}. 
On the other hand, when the intrinsic structure is considered, one can view it as an unconstrained optimization problem on the Stiefel manifolds. Then various kinds of unconstrained optimization algorithms can be extended to the corresponding Riemannian versions, \textit{e.g.}, \emph{steepest descent} methods \cite{abrudan2008steepest}, \emph{Barzilai-Borwein} methods \cite{Jiang2015A,Wen2013A} and \emph{trust region} methods \cite{Absil2007Trust,Ishteva2011Best}. 
For more details, one can refer to the monographs on Riemannian optimization \cite{absil2009optimization,Boumal2020An}. 
In particular, based on the relationship between homogeneous polynomial optimization and tensor computation, many Riemannian algorithms were used to solve the tensor decomposition problems, \textit{e.g.}, the \emph{geometric Newton} \cite{Elden2009A} and \emph{L-BFGS} \cite{Hans2018Nonlinearly} algorithms for Tucker decomposition \cite{kolda2009tensor}, and the polar decomposition based algorithms for low rank orthogonal approximation \cite{chen2009tensor,hu2022linear,Yang2020The,li2019polar}.
\fin

\subsection{Jacobi-type algorithms} 

If $r=n$ in problem \eqref{definition-f-new}, or $r_p=n_p$ for $1\leq p\leq d$ in problem \eqref{definition-f-prob-1}, then these problems have orthogonal constraints on unitary groups, or on the product of unitary groups.
In this case, besides the above optimization methods, one important approach is to use the Jacobi-type rotations \cite{Cardoso93:JADE,Como94:sp,Lathauwer96:simultaneous,Lathauwer01:ICA,martin2008jacobi,li2020convergence}.
From the computational point of view, a key step of Jacobi-type methods is to solve the subproblem.
For example, for the order-$3$ real symmetric tensor case of joint approximate symmetric tensor diagonalization, De Lathauwer, De Moor and Vandewalle \cite{Lathauwer01:ICA} showed that the solution of the subproblem has an SVD solution of a symmetric $2\times 2$ (real case) or $3\times 3$ (complex case) matrix. The order-$3$ real tensor (not necessarily symmetric) generalization of the Jacobi SVD algorithm \cite{GoluV96:jhu} also uses a similar method, while the subproblem is to maximize the trace of a tensor \cite{martin2008jacobi}. In this criterion, the subproblem has a straightforward explicit solution. If we set up the criterion as maximizing the sum of squares of the diagonal elements, it does not have a straightforward explicit solution, but still has an explicit solution involving the SVD. In the joint SVD problem for $2\times 2$ real matrices, it was shown in \cite{Pesquet2001Joint} that maximizing the sum of squares of the diagonal elements is equivalent to maximizing the trace. However, these two subproblem criteria are not equivalent in higher order tensor cases. 

For the best low multilinear rank
approximation of order-$d$ real symmetric tensors, Ishteva, Absil and Van Dooren \cite{IshtAV13:simax} proposed a \emph{gradient-based Jacobi-type} (Jacobi-G) algorithm, and proved the weak convergence\footnote{Every accumulation point is a stationary point.}.
Then, for the joint approximate symmetric tensor diagonalization on an orthogonal group with real symmetric matrices or order-$3$ tensors, Li, Usevich and Comon \cite{LUC2017globally} established the global convergence\footnote{For any starting point, the iterations converge as a whole sequence.} of Jacobi-G algorithm without any condition. For the joint approximate symmetric tensor diagonalization on a unitary group with complex matrices or tensors, Usevich, Li and Comon \cite{ULC2019} similarly formulated the Jacobi-G algorithm, established its global convergence under certain conditions, and estimated the local linear convergence rate.
For the joint approximate symmetric tensor diagonalization on a Stiefel manifold with real symmetric matrices or tensors, Li, Usevich and Comon \cite{li2019jacobi} formulated it as an equivalent optimization problem on an orthogonal group, and then used a Jacobi-type algorithm to solve it. 
To our knowledge, the global convergence of the Jacobi-type methods \cite{martin2008jacobi} for joint approximate (nonsymmetric) tensor diagonalization has not yet been established in the literature.
 
\subsection{Contributions}
In this paper, for the objective function $g$ in \eqref{definition-f-new}, we prove that the elementary function $\hij{i}{j}{\matr{U}}$ (see \eqref{eq-func-h-new}) can always be represented as the sum of a finite number of quadratic forms (see \cref{le:sub-problem-general}).  
If $\kappa=1$, this elementary function can always be represented as a quadratic form (see \cref{lem:coro-quadra-01}). 
For the Jacobi-G algorithm (\cref{Jacobi-G-general}) on a unitary group proposed in \cite{ULC2019}, we prove its global convergence under any one of two new conditions (see (C2) and (C3) in \cref{subsec:converg_analy}), when the subproblem can always be represented as a quadratic form. 
This global convergence result applies to the \emph{joint approximate symmetric tensor diagonalization} (JATD-S) (see \cref{ex:JATM-S-U}). 
For the Jacobi-MG algorithm (\cref{al:Jacobi-G-general}) on the product of unitary groups, we propose this algorithm, and prove that it is well defined.
Similar as for the Jacobi-G algorithm, we establish the global convergence of Jacobi-MG algorithm under any one of three conditions (see (C1), (C2) and (C3) in \cref{subsec:converg_analy}). 
This algorithm and its global convergence result apply to the \emph{joint approximate tensor diagonalization} (JATD) (see \cref{ex:nonsymmetric-diag-compression}). 
We also propose two general proximal variants: Jacobi-GP (\cref{Jacobi-GV-general}) and Jacobi-MGP (\cref{Jacobi-GV-general-gg}), and establish their weak and global convergence without any further condition. 
This paper is based on the complex Stiefel manifolds and complex tensors.
In fact, for the real cases of HPOSM problem \eqref{definition-f-new} and HPOSM-P problem \eqref{definition-f-prob-1}, all the algorithms and convergence results of this paper also apply. 

\subsection{Organization}

The paper is organized as follows. 
In \Cref{sec:geometri_unitary}, we recall the basics of geometries on the Stiefel manifolds and unitary groups, as well as the \L{}ojasiewicz gradient inequality.
 In \Cref{sec:back_summary}, we present the tensor approximation examples we study, and some geometric equivalence results. 
In \cref{sec:element}, we give the representation of the elementary functions under Jacobi rotations by using quadratic forms, as well as the Riemannian gradient. In \cref{sec:Jacobi-G algorithm}, we recall the Jacobi-G algorithm on the unitary group, and prove its global convergence under two different conditions. In \cref{sec:Jacobi-MG algorithm}, we develop a gradient-based multi-block Jacobi-type algorithm on the product of unitary groups, and prove its global convergence under three different conditions. In \cref{sec:Jacobi-GP algorithm}, we propose two proximal variants of the Jacobi-type algorithm, and establish their global convergence without any further condition. Some numerical experiments are shown in \cref{sec:experiments}. 
\cref{sec:conclusions} concludes this paper with some final remarks and possible future work.

\section{Geometries on the Stiefel manifold and unitary group}\label{sec:geometri_unitary}

\subsection{Wirtinger calculus}

For $\matr{X} \in \CC^{n\times r}$, we write $\matr{X} =\matr{X}^{\R} + \ui \matr{X}^{\I}$ for the real and imaginary parts.
For $\matr{X},\matr{Y}\in \CC^{n\times r}$, we introduce the following real-valued inner product
\begin{equation}\label{eq:RInner}
\RInner{\matr{X}}{\matr{Y}} \eqdef
\langle\matr{X}^{\R},\matr{Y}^{\R} \rangle + \langle\matr{X}^{\I},\matr{Y}^{\I} \rangle =
\Re \left({\trace{\matr{X}^{\H}\matr{Y}}}\right),
\end{equation}
which makes $\CC^{n\times r}$ a real Euclidean space of dimension $2nr$.
Let $g:\text{St}(r,n,\CC)\rightarrow\RR$ be a differentiable function and $\matr{X}\in\text{St}(r,n,\CC)$.
We denote by $\frac{\partial g}{\partial\matr{X}^\R},\frac{\partial g}{\partial\matr{X}^\I} \in \RR^{n\times r}$ the matrix Euclidean derivatives of $g$ with respect to real and imaginary parts of $\matr{X}$.
The \emph{Wirtinger derivatives} \cite{abrudan2008steepest,brandwood1983complex,krantz2001function} are defined as
\begin{equation*}
\frac{\partial g}{\partial\matr{X}^{*}} \eqdef \frac{1}{2}\left(\frac{\partial g}{\partial\matr{X}^\R}+ \ui\frac{\partial g}{\partial\matr{X}^\I}\right), \quad
\frac{\partial g}{\partial\matr{X}} \eqdef \frac{1}{2}\left(\frac{\partial g}{\partial\matr{X}^\R}- \ui\frac{\partial g}{\partial\matr{X}^\I}\right).
\end{equation*}
Then the Euclidean gradient of $g$ with respect to the inner product \eqref{eq:RInner} becomes
\begin{equation}\label{eq:matr_Euci_grad}
\nabla g(\matr{X}) = \frac{\partial g}{\partial\matr{X}^\R}+ \ui\frac{\partial g}{\partial\matr{X}^\I} = 2 \frac{\partial g}{\partial\matr{X}^{*}}.
\end{equation}
For real matrices $\matr{X},\matr{Y}\in\RR^{n\times r}$, we see that \eqref{eq:RInner} becomes the standard inner product, and \eqref{eq:matr_Euci_grad} becomes the standard Euclidean gradient.

\subsection{Riemannian gradient on \texorpdfstring {$\text{St}(r,n,\CC)$}{}}\label{subsec:riem-gradient}
We denote
$\skeww{\matr{P}}\eqdef\frac{1}{2}(\matr{P}-\matr{P}^{\H})$ 
for a complex matrix $\matr{P}\in\CC^{r\times r}$. 
Let ${\rm\bf T}_{\matr{X}} \text{St}(r,n,\CC)$ be the \emph{tangent space} to $\text{St}(r,n,\CC)$ at a point $\matr{X}$.
By \cite[Definition 6]{Manton2001Modified}, we know that
\begin{align}
{\rm\bf T}_{\matr{X}} \text{St}(r,n,\CC)
= \{\matr{Z}\in\CC^{n\times r}: \matr{Z}=\matr{X}\matr{A}+\matr{X}_{\perp}\matr{B},
\matr{A}\in\CC^{r\times r}, \matr{A}^{\H}+\matr{A}=0, \matr{B}\in\CC^{(n-r)\times r}\},\label{eq:tanget_spac}
\end{align}
which is a $(2nr-r^2)$-dimensional vector space.
The orthogonal projection of $\xi\in\CC^{n\times r}$ to \eqref{eq:tanget_spac} is
\begin{equation}\label{eq:proj_Stiefel}
{\rm Proj}_{\matr{X}} \xi = (\matr{I}_{n}-\matr{X}\matr{X}^{\H})\xi + \matr{X}\skeww{\matr{X}^{\H}\xi}.
\end{equation}
We denote ${\rm Proj}^{\bot}_{\matr{X}} \xi =\xi - {\rm Proj}_{\matr{X}} \xi $, which is in fact the orthogonal projection of $\xi $ to the \emph{normal space}\footnote{This is the orthogonal complement of the tangent space; see \cite[Section 3.6.1]{absil2009optimization} for more details.}.
Note that $\text{St}(r,n,\CC)$ is an embedded submanifold of the Euclidean space $\CC^{n\times r}$.
Let $g$ be the function in \eqref{definition-f-new}. 
From \eqref{eq:proj_Stiefel}, we have the \emph{Riemannian gradient}\footnote{See \cite[Equations (3.31) and (3.37)]{absil2009optimization} for a detailed definition.} of $g$ at $\matr{X}$ as:
\begin{align}\label{eq:Rie_grad}
\ProjGrad{g}{\matr{X}} = {\rm Proj}_{\matr{X}} \nabla g(\matr{X})
= (\matr{I}_{n}-\matr{X}\matr{X}^{\H}) \nabla g(\matr{X}) + \matr{X}\skeww{\matr{X}^{\H} \nabla g(\matr{X})}.
\end{align}
For the objective function $f$ in \eqref{definition-f-prob-1},
which is defined on the product of Stiefel manifolds, its Riemannian gradient can be computed as
$\ProjGrad{f}{\vect{\omega}} = (\ProjGrad{f_{(1)}}{\matr{X}^{(1)}},\ProjGrad{f_{(2)}}{\matr{X}^{(2)}},\cdots,\ProjGrad{f_{(d)}}{\matr{X}^{(d)}}).$ 

\subsection{Riemannian gradient on \texorpdfstring {$\UN{n}$}{}}
In \cref{subsec:riem-gradient}, if $r=n$, it specializes as the unitary group $\UN{n}$. In particular, the tangent space \cite[Section 3.5.7]{absil2009optimization} to $\UN{n}$ at a point $\matr{U}$ is 
\begin{eqnarray*}
\begin{aligned}
\mathbf{T}_{\matr{U}}\UN{n} = \{ \matr{Z}\in\CC^{n\times n} : \matr{Z}^{\H} \matr{U} + \matr{U}^{\H} \matr{Z} = \matr{0}\} 
= \{ \matr{Z}\in\CC^{n\times n} : \matr{Z} = \matr{U} \matr{P}, \,\, \matr{P} + \matr{P}^{\H} = \matr{0}, \,\, \matr{P}\in\CC^{n\times n}  \}.
\end{aligned}
\end{eqnarray*}
%
Let $\tilde{g}:\UN{n}\rightarrow\RR$ be a differentiable function and $\matr{U}\in\UN{n}$.
The equation \eqref{eq:Rie_grad} specializes as
\begin{eqnarray}\label{Rie-gradient-general}
\ProjGrad{\tilde{g}}{\matr{U}}=\matr{U} \skeww{\matr{U}^{\H} \nabla \tilde{g}(\matr{U})}=\matr{U} \matr{\Lambda}(\matr{U}),
\end{eqnarray}
 where $\matr{\Lambda}(\matr{U})\eqdef\skeww{\matr{U}^{\H} \nabla \tilde{g}(\matr{U})}$ is a skew-Hermitian matrix. 
For the objective function $\tilde{f}$ in \eqref{definition-f-ext-2}, 
its Riemannian gradient is
$\ProjGrad{\tilde{f}}{\vect{\upsilon}} = (\ProjGrad{\tilde{f}_{(1)}}{\matr{U}^{(1)}},\ProjGrad{\tilde{f}_{(2)}}{\matr{U}^{(2)}},\cdots,\ProjGrad{\tilde{f}_{(d)}}{\matr{U}^{(d)}}).$ 

\subsection{\L{}ojasiewicz gradient inequality}
We now recall some results about the \L{}ojasiewicz gradient inequality \cite{AbsMA05:sjo,loja1965ensembles,lojasiewicz1993geometrie,Usch15:pjo}, which will help us to establish the global convergence of Jacobi-type algorithms. 

\begin{definition} [{\cite[Definition 2.1]{SU15:pro}}]\label{def:Lojasiewicz}
Let $\mathcal{M} \subseteq \RR^n$ be a Riemannian submanifold and $f: \mathcal{M} \to \RR$ be a differentiable function. The function $f: \mathcal{M} \to \RR$ is said to satisfy a \emph{\L{}ojasiewicz gradient inequality} at $\vect{x} \in \mathcal{M}$, if there exist $\varpi>0$, $\zeta\in (0,\frac{1}{2}]$ and a neighborhood $\mathcal{U}$ in $\mathcal{M}$ of $\vect{x}$ such that
\begin{equation}\label{eq:Lojasiewicz}
|{f}(\vect{y})-{f}(\vect{x})|^{1-\zeta}\leq \varpi \|\ProjGrad{f}{\vect{y}}\|,
\end{equation}
for all $\vect{y}\in\mathcal{U}$.
\end{definition}
\begin{lemma}[{\cite[Proposition 2.2]{SU15:pro}}]\label{lemma-SU15}
Let $\mathcal{M}\subseteq\RR^n$ be an analytic submanifold\footnote{See {\cite[Definition 2.7.1]{krantz2002primer}} or \cite[Definition 5.1]{LUC2017globally} for a definition of an analytic submanifold.} and $f: \mathcal{M} \to \RR$ be a real analytic function. Then for any $\vect{x}\in \mathcal{M}$, $f$ satisfies a \L{}ojasiewicz gradient inequality \eqref{eq:Lojasiewicz} in the $\delta$-neighborhood of $\vect{x}$, for some\footnote{The values of $\delta, \varpi,\zeta$ depend on the specific point in question.} $\delta, \varpi>0$ and $\zeta\in (0,\frac{1}{2}]$.	
\end{lemma}
\begin{theorem}[{\cite[Theorem 2.3]{SU15:pro}}]\label{theorem-SU15}
Let $\mathcal{M}\subseteq\RR^n$ be an analytic submanifold and $\{\vect{x}_k\}_{k\geq 1}\subseteq\mathcal{M}$. Suppose that $f$ is a real analytic function, and for large enough $k$,\\
(i) there exists $\varsigma>0$ such that
\begin{equation*}\label{eq:sufficient_descent}
|{f}(\vect{x}_{k+1})-{f}(\vect{x}_k)|\geq \varsigma \|\ProjGrad{f}{\vect{x}_k}\|\|\vect{x}_{k+1}-\vect{x}_{k}\|
\end{equation*}
holds;\\
(ii) $\ProjGrad{f}{\vect{x}_k}=0$ implies that $\vect{x}_{k+1}=\vect{x}_{k}$.\\
Then any accumulation point $\vect{x}_*$ of $\{\vect{x}_k\}_{k\geq 1}$ must be the only limit point. 
\end{theorem}

\section{Tensor approximation examples and geometric results}\label{sec:back_summary}

\subsection{Examples in tensor approximation}\label{sec:problem-formulation-example}
We now present two important examples in tensor approximation, which can be written as the form \eqref{definition-f-new} or \eqref{definition-f-prob-1}. 

\begin{example}\label{ex:JATD-S}
Let $\{\tens{A}^{(\ell)}\}_{1\leq\ell\leq L}\subseteq\CC^{n\times n\times\cdots\times n}$ be a set of order-$d$ tensors.
Let $\alpha_\ell \in \RR^{+}$ for $1\leq\ell\leq L$.
Let $1\leq r\leq n$. 
The \emph{joint approximate symmetric tensor diagonalization} (JATD-S) \cite{chen2009tensor,LUC2017globally,li2019jacobi,li2019polar} is
\begin{align}\label{cost-fn-general-1-s}
\max\ g(\matr{X}) = \sum\limits_{\ell=1}^{L} \alpha_\ell\|\diag{\tens{W}^{(\ell)}}\|^2,\ \ \matr{X}\in\text{St}(r,n,\CC),
\end{align}
where $\tens{W}^{(\ell)}=\tens{A}^{(\ell)} \contr{1} \matr{X}^\dagger\contr{2}\cdots \contr{d} \matr{X}^\dagger$ for $1\leq \ell\leq L$, and  $(\cdot)^{\dagger}$ stands for $(\cdot)^{\H}$ or $(\cdot)^{\T}$. 
Then we have that
\begin{align}\label{cost-fn-general-1-s-3}
g(\matr{X}) = \sum\limits_{q=1}^{r}\sum_{\ell=1}^{L} \alpha_\ell|\tens{A}^{(\ell)}\contr{1}\vect{x}_{q}^\dagger\contr{2}\cdots\contr{d}\vect{x}_{q}^\dagger|^2= \sum\limits_{q=1}^{r}\matr{\Gamma}\left(\left(\sum_{\ell=1}^{L} \alpha_\ell{\tens{A}^{(\ell)}}^{*}{\otimes\tens{A}^{(\ell)}}\right)^{\sharp},\vect{x}_{q},d\right), 
\end{align}
and thus it is of the form \eqref{definition-f-new}
with $\kappa=d$. 
Here, for a complex matrix $\tens{B}$, we denote that
\begin{equation}\label{eq:def_sharp}
(\tens{B})^{\sharp}\eqdef\bigg\{
\begin{aligned}
\tens{B}^{*},\ \ & \textrm{if}\ (\cdot)^{\dagger}=(\cdot)^{\H}, \\
\tens{B},\ \ & \textrm{if}\ (\cdot)^{\dagger}=(\cdot)^{\T}.
\end{aligned}
\end{equation}
\end{example}
\begin{example}\label{ex:JATD}
Let $\{\tens{A}^{(\ell)}\}_{1\leq\ell\leq L}\subseteq\CC^{n_1\times n_2\times\cdots n_d}$ be a set of complex tensors.
Let $\alpha_\ell \in \RR^{+}$ for $1\leq\ell\leq L$.
Let $1\leq r\leq n_p$ for $1\leq p\leq d$.
The \emph{joint approximate tensor diagonalization} (JATD) \cite{chen2009tensor,li2019polar} problem is
\begin{align}\label{cost-fn-general-1}
\max\ f(\vect{\omega}) = \sum\limits_{\ell=1}^{L} \alpha_\ell\|\diag{\tens{W}^{(\ell)}}\|^2,\ \ \matr{X}^{(p)}\in\text{St}(r,n_p,\CC),
\end{align}
where $\tens{W}^{(\ell)}=\tens{A}^{(\ell)} \contr{1} (\matr{X}^{(1)})^{\dagger}\contr{2}\cdots \contr{d} (\matr{X}^{(d)})^{\dagger}$ for $1\leq \ell\leq L$, and $(\cdot)^{\dagger}$ stands for $(\cdot)^{\H}$ or $(\cdot)^{\T}$. 
If we fix $1\leq p\leq d$, and denote
\begin{align*}
\tens{V}^{(\ell)}&= \tens{A}^{(\ell)} \contr{1} (\matr{X}^{(1)})^{\dagger}\contr{2}\cdots \contr{p-1} (\matr{X}^{(p-1)})^{\dagger}\contr{p+1} (\matr{X}^{(p+1)})^{\dagger}\cdots\contr{d} (\matr{X}^{(d)})^{\dagger},\\
\vect{v}^{(\ell)}_{q}&=\tens{A}^{(\ell)}\contr{1}(\vect{x}_{q}^{(1)})^\dagger\contr{2}\cdots\contr{p-1}(\vect{x}_{q}^{(p-1)})^\dagger\contr{p+1}(\vect{x}_{q}^{(p+1)})^\dagger\cdots\contr{d}(\vect{x}_{q}^{(d)})^\dagger,
\end{align*}
then the restricted function $f_{(p)}$ in \eqref{definition-h} can be represented as
\begin{align}\label{cost-fn-general-1-5}
f_{(p)}(\matr{X}) = \sum_{q=1}^{r}\sum\limits_{\ell=1}^{L} \alpha_\ell|\vect{x}_{q}^\dagger\vect{v}^{(\ell)}_{q}|^2
= \sum_{q=1}^{r}\left(\sum\limits_{\ell=1}^{L} \alpha_\ell(\vect{v}^{(\ell)}_{q})^{*}(\vect{v}^{(\ell)}_{q})^\T\right)^{\sharp} \contr{1}\vect{x}_{q}^\H\contr{2}\vect{x}_{q}^\T,
\end{align}
where $(\cdot)^{\sharp}$ is as in \eqref{eq:def_sharp}, 
and thus it is of the form \eqref{definition-f-new} with $\kappa=1$.
\end{example}

\subsection{Problem reformulations on \texorpdfstring {$\UN{n}$}{} and \texorpdfstring {$\Upsilon$}{}}

\subsubsection{For problem \eqref{definition-f-new} on a Stiefel manifold}
Let
\begin{equation}\label{definition-f-g}
g:\ \text{St}(r,n,\CC) \longrightarrow \RR
\end{equation}
be an abstract function, which includes the objective function \eqref{definition-f-new} as a special case. 
 This paper aims to use Jacobi-type elementary rotations to optimize the abstract function $g$ in \eqref{definition-f-g}. 
To do this, we first need to reformulate problem \eqref{definition-f-g} from the Stiefel manifold $\text{St}(r,n,\CC)$ to the \emph{unitary group} $\UN{n}$.
Let $\sigma:\UN{n}\rightarrow\text{St}(r,n,\CC)$ be the natural projection map which only keeps the first $r$ columns.
We define 
\begin{equation}\label{definition-f-ext}
\tilde{g}:\ \UN{n} \longrightarrow \RR,\,\,\,\, \vect{\matr{U}}\mapsto g(\sigma(\vect{\matr{U}})).
\end{equation}
 An equivalence relationship between the stationary points of functions $g$ in \eqref{definition-f-g} and $\tilde{g}$ in \eqref{definition-f-ext} will be shown in \cref{subsec:geometr-equiva} (see \cref{theorem-equiva-stationary}). 
\fin

\begin{example}\label{ex:JATM-S-U}
Let $g$ be the homogeneous polynomial in \eqref{definition-f-new}. 
Denote $\matr{U}=(\vect{u}_1,\vect{u}_2,\cdots,\vect{u}_n)\in\UN{n}$. 
Then the corresponding $\tilde{g}$ can be represented as 

\begin{align}\label{definition-f-new-32}
\tilde{g}(\matr{U}) = \sum_{q=1}^{r}
\matr{\Gamma}(\tens{B}^{(q)},\vect{u}_{q},\kappa).
\end{align}
\fin
If $g$ is the example function in \eqref{cost-fn-general-1-s}, then we have 
\begin{align}
\tilde{g}(\matr{U}) = \sum\limits_{\ell=1}^{L} \alpha_\ell\|\diag{(\tens{W}^{(\ell)})^{(r,r,\cdots, r)}} \|^2,
\end{align}
where $\tens{W}^{(\ell)}=\tens{A}^{(\ell)} \contr{1} \matr{U}^\dagger\contr{2}\cdots \contr{d} \matr{U}^\dagger$ for $1\leq \ell\leq L.$
\end{example}

\subsubsection{For problem \eqref{definition-f-prob-1} on the product of Stiefel manifolds}
Let
\begin{equation}\label{definition-f-g-23}
f:\ \Omega \longrightarrow \RR
\end{equation}
be an abstract function, which includes the objective function \eqref{definition-f-prob-1} as a special case.
We denote that 
$\Upsilon\eqdef\UN{n_1}\times\UN{n_2}\times\cdots\times\UN{n_d}$, 
and denote by $\rho:\Upsilon\rightarrow\Omega$ the natural projection map, which only keeps the first $r_p$ columns for the $p$-th block variable. 
Denote $\vect{\upsilon}=(\matr{U}^{(1)}, \matr{U}^{(2)},\cdots,\matr{U}^{(d)})\in\Upsilon$.
For problem \eqref{definition-f-g-23}, we define
\begin{equation}\label{definition-f-ext-2}
\tilde{f}:\ \Upsilon \longrightarrow \RR,\ \vect{\upsilon}\mapsto f(\rho(\vect{\upsilon})).
\end{equation}
 Similarly, in \cref{subsec:geometr-equiva}, we will prove an equivalence relationship between the stationary points of functions $f$ in \eqref{definition-f-g-23} and $\tilde{f}$ in \eqref{definition-f-ext-2} (see \cref{theorem-equiva-stationary-2}).

\begin{example}\label{ex:nonsymmetric-diag-compression}
Let $f$ be the example function in \eqref{cost-fn-general-1}. Then we have 
\begin{align}\label{eq:f-tilde-diag}
\tilde{f}(\vect{\upsilon}) = \sum\limits_{\ell=1}^{L} \alpha_\ell\| \diag{(\tens{W}^{(\ell)})^{(r,r,\cdots,r)}} \|^2,\ \ \matr{U}^{(p)}\in\UN{n_p},\ 1\leq p\leq d,
\end{align}
 where $\tens{W}^{(\ell)}=\tens{A}^{(\ell)} \contr{1} (\matr{U}^{(1)})^{\dagger}\contr{2}\cdots \contr{d} (\matr{U}^{(d)})^{\dagger}$ for $1\leq \ell\leq L$. 
\end{example}

\subsection{Geometric equivalence}\label{subsec:geometr-equiva}
Let $g$ be as in \eqref{definition-f-g}, and $\tilde{g}$ be as in \eqref{definition-f-ext}. 
Let $\mathbb{S}\subseteq\CC$ be the set of unimodular complex numbers.
We say that $g$ is \emph{scale invariant}, if
$g(\matr{X}) = g(\matr{X}\matr{R})$ for all $\matr{X}\in\text{St}(r,n,\CC)$ and
{\small
\begin{equation*}\label{eq:invariance_scaling-0}
\matr{R}\in\set{DU}_{r} \eqdef\left\{
\begin{bmatrix}z_1 & & 0 \\ & \ddots & \\ 0 & & z_r \end{bmatrix}: z_p\in\mathbb{S},\ 1\leq p\leq r\right\}\subseteq\UN{r}.
\end{equation*}}
Let $f$ be as in \eqref{definition-f-g-23} and $\tilde{f}$ be as in \eqref{definition-f-ext-2}. We say that $f$ is \emph{scale invariant}, if for all $1\leq p\leq d$, the restricted function $f_{(p)}$ in \eqref{definition-h} is always scale invariant.
 We now show some interesting relationships between the geometric properties of $g$ and $\tilde{g}$, as well as that between $f$ and $\tilde{f}$.
\begin{theorem}\label{theorem-equiva-invar}
 Function $g$ is scale invariant if and only if $\tilde{g}$ is scale invariant.
\end{theorem}

\begin{proof}
If $g$ is scale invariant, $\matr{U}\in\UN{n}$ and $\matr{R}\in\set{DU}_{n}$,
then
\begin{align*}
\tilde{g}(\matr{U}\matr{R}) = g(\sigma(\vect{\matr{U}\matr{R}})) = g(\sigma(\vect{\matr{U}})\matr{R}^{(r,r)})=g(\sigma(\vect{\matr{U}})) = \tilde{g}(\matr{U}),
\end{align*}
and thus $\tilde{g}$ is also scale invariant. Conversely, if $\tilde{g}$ is scale invariant, $\matr{X}\in\text{St}(r,n,\CC)$ and $\matr{R}\in\set{DU}_{r}$,
then there exists $\matr{U}\in\UN{n}$ satisfying $\sigma(\matr{U})=\matr{X}$, and $\matr{R}_{0}\in\set{DU}_{n}$ satisfying $\matr{R}_{0}^{(r,r)}=\matr{R}$. Then
\begin{align*}
g(\matr{X}\matr{R}) = g(\sigma(\matr{U})\matr{R}) = g(\sigma(\matr{U}\matr{R}_{0})) = \tilde{g}(\matr{U}\matr{R}_{0}) = \tilde{g}(\matr{U})=g(\matr{X}).
\end{align*}
The proof is complete.
\end{proof}

The following result is a direct consequence of \cref{theorem-equiva-invar}.

\begin{corollary}
 Function $f$ is scale invariant if and only if $\tilde{f}$ is scale invariant.
\end{corollary}

\begin{remark}
It is clear that the functions $g$ in \eqref{definition-f-new} and $f$ in \eqref{definition-f-prob-1} are both scale invariant.
Therefore, both the example functions in \cref{sec:problem-formulation-example} are scale invariant.
\end{remark}

\begin{theorem}\label{theorem-equiva-stationary}
Let $\matr{U}_*\in \UN{n}$ and $\matr{X}_*\in\text{St}(r,n,\CC)$ satisfy $\sigma(\matr{U}_*)=\matr{X}_*$.
Then $\ProjGrad{\tilde{g}}{\matr{U}_{*}}=\matr{0}$ if and only if $\ProjGrad{g}{\matr{X}_{*}}=\matr{0}$.
\end{theorem}

\begin{proof}
 Denote $\matr{U}_*= [\matr{X}_* \ \matr{Y}_*] \in \UN{n}$ with $\matr{X}_*\in\text{St}(r,n,\CC)$. Then we have  
\begin{eqnarray}\label{eq:Ug}
\begin{aligned}
\matr{U}_*^{\H} \nabla \tilde{g}(\matr{U}_*) &= \left[ \begin{array}{c}
\matr{X}_*^{\H} \\
\matr{Y}_*^{\H} \end{array} \right] \left[\partial \tilde{g}(\matr{X}_*),\,\, \partial \tilde{g}(\matr{Y}_*) \right]
= \begin{bmatrix} \matr{X}_*^{\H} \partial \tilde{g}(\matr{X}_*) & \matr{X}_*^{\H} \partial \tilde{g}(\matr{Y}_*) \\ \matr{Y}_*^{\H} \partial \tilde{g}(\matr{X}_*) & \matr{Y}_*^{\H} \partial \tilde{g}(\matr{Y}_*) \end{bmatrix}.
\end{aligned}
\end{eqnarray}
If $\ProjGrad{g}{\matr{X}_{*}}=\matr{0}$, equation \eqref{eq:Rie_grad} tells us that  
\begin{eqnarray}\label{eq:St-Rie-g}
\left( \matr{I}_{n}-\matr{X}_*\matr{X}_*^{\H} \right)\nabla g(\matr{X}_*) + \matr{X}_*\skeww{\matr{X}_*^{\H}\nabla g(\matr{X}_*)}
=\matr{0}.
\end{eqnarray} 
\fin Multiplying both sides of \eqref{eq:St-Rie-g} by $\matr{X}_*^{\H}$ to the left,
we get $\skeww{\matr{X}_*^{\H}\nabla g(\matr{X}_*)}=\matr{0}$. 
It follows from \eqref{eq:St-Rie-g} that $\left( \matr{I}_{n}-\matr{X}_*\matr{X}_*^{\H} \right)\nabla g(\matr{X}_*)=\matr{0}$, and thus 
\begin{equation*}\label{eq:XY0}
\matr{Y}_*\matr{Y}_*^{\H} \partial \tilde{g}(\matr{X}_*)= \left( \matr{I}_{n}-\matr{X}_*\matr{X}_*^{\H} \right) \nabla g(\matr{X}_*) = \matr{0},
\end{equation*}
which implies that  $\matr{Y}_*^{\H}\partial \tilde{g}(\matr{X}_*)=\matr{0}$. 
Note that $\partial \tilde{g}(\matr{Y}_*)=\matr{0}$ and we have proved $\skeww{\matr{X}_*^{\H}\nabla g(\matr{X}_*)}=\matr{0}$.
By \eqref{eq:Ug}, 
we see that $\skeww{\matr{U}_*^{\H}\nabla \tilde{g}(\matr{U}_*)}=\matr{0}$, and then 
$\ProjGrad{\tilde{g}}{\matr{U}_*}=\matr{U}_*\skeww{\matr{U}_*^{\H}\nabla \tilde{g}(\matr{U}_*)}
=\matr{0}.$  
Conversely, if $\ProjGrad{\tilde{g}}{\matr{U}_{*}}=\matr{0}$, by \eqref{eq:Ug}, we see that $\skeww{\matr{X}_*^{\H}\nabla g(\matr{X}_*)}=\matr{0}$ and $\matr{Y}_*^{\H}\partial \tilde{g}(\matr{X}_*)=\matr{0}$. It follows that $\left( \matr{I}_{n}-\matr{X}_*\matr{X}_*^{\H} \right)\nabla g(\matr{X}_*)=\matr{0}$, and thus $\ProjGrad{g}{\matr{X}_*} = \matr{0}$. 
The proof is complete. 
\end{proof}

Similarly, we also have the following result for $f$ in \eqref{definition-f-g-23} and $\tilde{f}$ in \eqref{definition-f-ext-2}. 

\begin{corollary}\label{theorem-equiva-stationary-2}
Let $\vect{\omega}\in\Omega$ and $\vect{\upsilon}\in\Upsilon$ satisfy $\rho(\vect{\upsilon})=\vect{\omega}$.
Then $\ProjGrad{\tilde{f}}{\vect{\upsilon}}=\matr{0}$ if and only if $\ProjGrad{f}{\vect{\omega}}=\matr{0}$.
\end{corollary}

\begin{remark}
In \cite[Proposition 2.8]{li2019jacobi}, a similar result as \cref{theorem-equiva-stationary} was proved for low rank orthogonal approximation problem of real symmetric tensors.
A more general result of \cref{theorem-equiva-stationary} can also be seen in \cite[Proposition 9.3.8]{Boumal2020An}.
\end{remark}

\section{Elementary functions}\label{sec:element}

\subsection{Elementary functions}\label{subsec:elementary}
Let $(i,j)$ be an index pair satisfying that $1\leq i<j\leq n$.
Let $\Gmat{i}{j}{\matr{\Psi}}$ be the \emph{plane transformation} defined as in \cite[Equation (2.8)]{ULC2019} for $\matr{\Psi} \in \UN{2}$.
%
Let $\tilde{g}$ be as in \eqref{definition-f-ext}.
Fix $\matr{U}\in\UN{n}$.
 Define the elementary function 
\begin{equation}\label{eq-func-h-new}
\begin{split}
\hij{i}{j}{\matr{U}}:\ \UN{2} \longrightarrow \RR, \ 
 {\matr{\Utwo}}\longmapsto \tilde{g}(\matr{U}\Gmat{i}{j}{\matr{\Utwo}}).
\end{split}
\end{equation} 
\fin
Note that the function $\tilde{g}$ in \eqref{definition-f-new-32} is scale invariant. We only need to consider
\begin{align}
\matr{\Psi}= \matr{\Utwo}(c,s_1,s_2)
&= \begin{bmatrix}
c & -s \\
s^{\ast} & c
\end{bmatrix}
= \begin{bmatrix}
c & -(s_1+\ui s_2) \\
s_1-\ui s_2 & c
\end{bmatrix}= \begin{bmatrix}
\cos\theta & -\sin\theta \ue^{\ui\phi} \\
\sin\theta \ue^{-\ui\phi} & \cos\theta
\end{bmatrix},
\label{unitary-para-2}
\end{align}
where $c\in \RR^{+}$, $s= s_1 + \ui s_2\in\CC$, $c^2+|s|^2 =1$.
For $1\leq \alpha, \beta\leq \kappa$, we denote
\begin{equation}\label{eq:omega_t}
\vect{z}_{\alpha,\beta} = (\cos \alpha\theta, -\sin \alpha\theta\cos\beta\phi, -\sin \alpha\theta\sin\beta\phi)^{\T} \in \mathbb{R}^3.
\end{equation} 
For $1\leq r\leq n$,
we define the pair index set $\mathbb{I} \eqdef\{(i,j):1\leq i<j \leq n\}$, and divide it into three different subsets
\begin{equation*}
\mathbb{I}_1=\{(i,j):1\leq i<j \leq r\}, \ \mathbb{I}_2=\{(i,j):1\leq i \leq r <j \leq n\}, \ \mathbb{I}_3=\{(i,j):r< i<j \leq n\}.
\end{equation*}
Let $\tilde{g}$ be as in \eqref{definition-f-new-32}.
For $1\leq q\leq r$,
we denote that
$\tens{W}^{(q)}\eqdef\matr{\Gamma}(\tens{B}^{(q)},\matr{U},\kappa)$ and $\tens{T}^{(q)} \eqdef (\tens{W}^{(q)})^{(i,j)}\in\CC^{2\times 2\times\cdots\times 2}$. 
Let $\Gmat{i}{j}{\matr{\Psi}}=(\vect{r}_1,\vect{r}_2,\cdots,\vect{r}_n)\in\UN{n}$, and
\begin{align}\label{eq:x-vect-y}
\vect{x} =
\begin{bmatrix}
\cos\theta \\
\sin\theta \ue^{-\ui\phi}
\end{bmatrix}, \ \
\vect{y} =
\begin{bmatrix}
-\sin\theta \ue^{\ui\phi} \\
\cos\theta
\end{bmatrix}.
\end{align}
Then the elementary function $\hij{i}{j}{\matr{U}}$ in \eqref{eq-func-h-new} can be further represented as 
\begin{itemize} 
\item if $(i,j)\in\mathbb{I}_1$, then
\begin{eqnarray*}
\begin{aligned}
\hij{i}{j}{\matr{U}}(\matr{\Utwo}) = \sum_{q=1}^{r}\matr{\Gamma}(\tens{W}^{(q)},\vect{r}_{q},\kappa)
= \matr{\Gamma}(\tens{T}^{(i)},\vect{x},\kappa) + \matr{\Gamma}(\tens{T}^{(j)},\vect{y},\kappa) +C_1,\label{eq:h-represen-u-01}
\end{aligned}
\end{eqnarray*}
where $C_1$ is a constant.
\item if $(i,j)\in\mathbb{I}_2$, then
$\hij{i}{j}{\matr{U}}(\matr{\Utwo}) =\matr{\Gamma}(\tens{T}^{(i)},\vect{x},\kappa)+C_2$, where $C_2$ is a constant. 
\item if $(i,j)\in\mathbb{I}_3$, then
$\hij{i}{j}{\matr{U}}(\matr{\Utwo}) = C_3$, where $C_3$ is a constant. 
\end{itemize}

\subsection{Representation using quadratic forms}\label{subsec:quad-form-complex}

In this subsection, we first prove that $\hij{i}{j}{\matr{U}}$ can always be represented as the sum of a finite number of quadratic forms. 
\begin{proposition}\label{le:sub-problem-general}
Let $\tilde{g}$ be as in \eqref{definition-f-new-32}. 
Let $\tens{T}$ be an order-$2\kappa$ 2-dimensional Hermitian tensor and $\kappa$-semisymmetric. 
Let $\vect{x}$ and $\vect{y}$ be as in \eqref{eq:x-vect-y}. 
 Then 
\begin{flalign}
&& \matr{\Gamma}(\tens{T},\vect{x},\kappa) &= \sum_{\alpha=1}^{\kappa}\sum_{\beta=1}^{\kappa} \vect{z}_{\alpha,\beta}^{\T}\matr{M}^{(\alpha,\beta)}\vect{z}_{\alpha,\beta} + C_1, & \label{eq:subproblem-general-x} \\
&& \matr{\Gamma}(\tens{T},\vect{y},\kappa) &= \sum_{\alpha=1}^{\kappa}\sum_{\beta=1}^{\kappa} \vect{z}_{\alpha,\beta}^{\T}(-1)^{\alpha}\matr{M}^{(\alpha,\beta)}\vect{z}_{\alpha,\beta} + C_2, & \label{eq:subproblem-general-y}
\end{flalign}
\fin
where the unit vectors $\vect{z}_{\alpha,\beta}$ are as in \eqref{eq:omega_t}, $\matr{M}^{(\alpha,\beta)}\in\RR^{3\times 3}$ are real symmetric matrices determined by $\tens{T}$, $C_1$ and $C_2$ are constants.
\end{proposition}

To prove \cref{le:sub-problem-general}, we need first to present some formulas of trigonometric functions, which can be easily proved by induction methods.

\begin{lemma}\label{le:equivalent-trigonometric-functions-even}
Let $k_1,k_2\geq 1$. Then we have 
%
$$\cos^{2k_1}\theta\sin^{2k_2}\theta = \sum_{t=0}^{k_1+k_2} c_t\cos 2t\theta \quad \text{and} \quad
\cos^{2k_2}\theta\sin^{2k_1}\theta = \sum_{t=0}^{k_1+k_2} (-1)^t c_t\cos 2t\theta,$$ 
%
where $\{c_t\}_{t=0}^{k_1+k_2}$ are constants satisfying that $\sum_{t=0}^{k_1+k_2} c_t = 0$ and $\sum_{t=0}^{k_1+k_2} (-1)^{t}c_t = 0$. 
\end{lemma}
\begin{lemma}\label{le:equivalent-trigonometric-functions}
Let $k_1,k_2\geq 1$. Then we have 
%
$$\cos^{2k_1+1}\theta\sin^{2k_2-1}\theta = \sum_{t=1}^{k_1+k_2} s_t\sin 2t\theta \quad \text{and} \quad 
\cos^{2k_2-1}\theta\sin^{2k_1+1}\theta = \sum_{t=1}^{k_1+k_2} (-1)^{t-1} s_t\sin 2t\theta,$$ 
%
where $\{s_t\}_{t=1}^{k_1+k_2}$ are constants. 
\end{lemma}

\begin{proof}[Proof of \cref{le:sub-problem-general}]
Denote $\tbinom{\kappa}{t} = \frac{\kappa!}{t!(\kappa-t)!}$ for $1\leq t\leq\kappa$, and $\tbinom{\kappa}{0}=\kappa$. 
We now first prove equation \eqref{eq:subproblem-general-x}.
Note that $\tens{T}$ is $\kappa$-semisymmetric.
We have that 
{\small\begin{align}
 \matr{\Gamma}(\tens{T},\vect{x},\kappa) &= \sum_{t_1,t_2=0}^{\kappa} \tbinom{\kappa}{t_1}\tbinom{\kappa}{t_2} \tenselem{T}_{\underbrace{1\cdots1}_{\kappa-t_1}\underbrace{2\cdots2}_{t_1}\underbrace{1\cdots1}_{\kappa-t_2}\underbrace{2\cdots2}_{t_2}}\cos^{2\kappa-(t_1+t_2)}\theta\sin^{t_1+t_2}\theta\ue^{(t_1-t_2)\ui\phi}\notag\\
& = \sum_{0\leq t_1\leq t_2\leq \kappa} \tbinom{\kappa}{t_1}\tbinom{\kappa}{t_2} \cos^{2\kappa-(t_1+t_2)}\theta\sin^{t_1+t_2}\theta\notag\\ &\quad\cdot\left(\tenselem{T}_{\underbrace{1\cdots1}_{\kappa-t_1}\underbrace{2\cdots2}_{t_1}\underbrace{1\cdots1}_{\kappa-t_2}\underbrace{2\cdots2}_{t_2}}\ue^{(t_1-t_2)\ui\phi} + \tenselem{T}_{\underbrace{1\cdots1}_{\kappa-t_2}\underbrace{2\cdots2}_{t_2}\underbrace{1\cdots1}_{\kappa-t_1}\underbrace{2\cdots2}_{t_1}}\ue^{(t_2-t_1)\ui\phi}\right)\notag\\
& = \sum_{0\leq t_1\leq t_2\leq \kappa} \tbinom{\kappa}{t_1}\tbinom{\kappa}{t_2} \cos^{2\kappa-(t_1+t_2)}\theta\sin^{t_1+t_2}\theta\notag\\ &\quad \cdot 2\left( \tenselem{T}_{\underbrace{1\cdots1}_{\kappa-t_1}\underbrace{2\cdots2}_{t_1}\underbrace{1\cdots1}_{\kappa-t_2}\underbrace{2\cdots2}_{t_2}}^{\R}\cos(t_2-t_1)\phi + \tenselem{T}_{\underbrace{1\cdots1}_{\kappa-t_1}\underbrace{2\cdots2}_{t_1}\underbrace{1\cdots1}_{\kappa-t_2}\underbrace{2\cdots2}_{t_2}}^{\I}\sin(t_2-t_1)\phi \right).\label{eq:c-x-1}
\end{align}}
If $t_2-t_1$ is even, by \cref{le:equivalent-trigonometric-functions-even}, we have %
\begin{align*}
\cos^{2\kappa-(t_1+t_2)}\theta\sin^{t_1+t_2}\theta\cos(t_2-t_1)\phi
&= \sum_{t=0}^{\kappa} c_t(1-2\sin^2 t\theta) \cos(t_2-t_1)\phi
= \sum_{t=0}^{\kappa} -2c_t\sin^2 t\theta \cos(t_2-t_1)\phi\\
&=\sum_{t=0}^{\kappa} -2c_t\sin^2 t\theta \left(\cos^2\frac{(t_2-t_1)\phi}{2}-\sin^2\frac{(t_2-t_1)\phi}{2}\right) ,\\
\cos^{2\kappa-(t_1+t_2)}\theta\sin^{t_1+t_2}\theta\sin(t_2-t_1)\phi
&= \sum_{t=0}^{\kappa} c_t(1-2\sin^2 t\theta) \sin(t_2-t_1)\phi= \sum_{t=0}^{\kappa} -2c_t\sin^2 t\theta \sin(t_2-t_1)\phi\\
&=\sum_{t=0}^{\kappa} -4c_t\sin^2 t\theta \cos\frac{(t_2-t_1)\phi}{2}\sin\frac{(t_2-t_1)\phi}{2}.
\end{align*}
If $t_2-t_1$ is odd, by \cref{le:equivalent-trigonometric-functions}, we have 
\begin{align*}
\cos^{2\kappa-(t_1+t_2)}\theta\sin^{t_1+t_2}\theta\cos(t_2-t_1)\phi
&= \sum_{t=1}^{\kappa} 2s_t\cos t\theta\sin t\theta \cos(t_2-t_1)\phi,\\
\cos^{2\kappa-(t_1+t_2)}\theta\sin^{t_1+t_2}\theta\sin(t_2-t_1)\phi
&= \sum_{t=1}^{\kappa} 2s_t\cos t\theta\sin t\theta \sin(t_2-t_1)\phi.
\end{align*}
It can be seen that every item in \eqref{eq:c-x-1} belongs to a quadratic form in \eqref{eq:subproblem-general-x}, and thus the equation \eqref{eq:subproblem-general-x} is proved. 
Similarly, we can prove equation \eqref{eq:subproblem-general-y} by \cref{le:equivalent-trigonometric-functions-even,le:equivalent-trigonometric-functions}. 
The proof is complete. 
\end{proof}
\fin

\begin{remark}\label{re:non-unique}
The matrices $\matr{M}^{(\alpha,\beta)}$ in \eqref{eq:subproblem-general-x} and \eqref{eq:subproblem-general-y} may be not unique. For example, we have $\cos 2\theta\sin 2\theta\cos 2\phi=2\cos\theta\sin\theta\cos^2\phi-\cos\theta\sin\theta$, since $\cos2\alpha\theta=2\cos^2\alpha\theta-1$ and $\sin2\alpha\theta=2\sin\alpha\theta\cos\alpha\theta$. 
\end{remark}

We now present two concrete examples of \cref{le:sub-problem-general}.

\begin{example}\label{lem:quadratic-gamma-01}
Let $\tens{T}\in\CC^{2\times 2}$ be a Hermitian matrix.
Let $\vect{x}$ and $\vect{y}$ be as in \eqref{eq:x-vect-y}.
Then \\
(i)
$\tens{T}\contr{1}\vect{x}^\H\contr{2}\vect{x}^\T= \vect{z}_{1,1}^\T \matr{M}^{(1,1)} \vect{z}_{1,1} + \tenselem{T}_{22}$, \quad \\
(ii) $\tens{T}\contr{1}\vect{y}^\H\contr{2}\vect{y}^\T= -\vect{z}_{1,1}^\T \matr{M}^{(1,1)} \vect{z}_{1,1} + \tenselem{T}_{11}$,\\
where  $$\matr{M}^{(1,1)}=
\begin{bmatrix}
\tenselem{T}_{11}-\tenselem{T}_{22} & -\tenselem{T}_{12}^{\R} & -\tenselem{T}_{12}^{\I}\\
-\tenselem{T}_{12}^{\R} & 0 & 0\\
-\tenselem{T}_{12}^{\I} & 0 & 0
\end{bmatrix}.$$ 
\end{example}

\begin{corollary}\label{lem:coro-quadra-01}
Let $\tilde{g}$ be as in \eqref{definition-f-new-32} with $\kappa=1$. Then the elementary function $\hij{i}{j}{\matr{U}}$ in \eqref{eq-func-h-new} can be expressed as a quadratic form if $(i,j)\in\mathbb{I}_1$ or $(i,j)\in\mathbb{I}_2$.
\end{corollary}

\begin{example}\label{lem:quadratic-gamma-02}
Let $\tens{T}\in\CC^{2\times 2\times 2\times 2}$ be an order-$4$ Hermitian tensor and $2$-semisymmetric.
Let $\vect{x}$ and $\vect{y}$ be as in \eqref{eq:x-vect-y}. 
Then \\
(i) $\matr{\Gamma}(\tens{T},\vect{x},2)= \vect{z}_{2,1}^\T \matr{M}^{(2,1)}
\vect{z}_{2,1} +\vect{z}_{1,1}^\T \matr{M}^{(1,1)}
\vect{z}_{1,1}+ \frac{1}{4}(3\tenselem{T}_{2222}-\tenselem{T}_{1111})$, \\
(ii) $\matr{\Gamma}(\tens{T},\vect{y},2)= \vect{z}_{2,1}^\T \matr{M}^{(2,1)}
\vect{z}_{2,1} -\vect{z}_{1,1}^\T \matr{M}^{(1,1)}
\vect{z}_{1,1}+\frac{1}{4}(3\tenselem{T}_{1111}-\tenselem{T}_{2222})$,\\
where
\begin{align*}
\matr{M}^{(2,1)} &= \frac{1}{4}
\begin{bmatrix}
\tenselem{T}_{1111}+\tenselem{T}_{2222} & 2\tenselem{T}_{1222}^{\R}-2\tenselem{T}_{1112}^{\R} & 2\tenselem{T}_{1222}^{\I}-2\tenselem{T}_{1112}^{\I}\\
2\tenselem{T}_{1222}^{\R}-2\tenselem{T}_{1112}^{\R} & 2\tenselem{T}_{1122}^{\R}+4\tenselem{T}_{1212} & 2 \tenselem{T}_{1122}^{\I}\\
2\tenselem{T}_{1222}^{\I}-2\tenselem{T}_{1112}^{\I} & 2 \tenselem{T}_{1122}^{\I} & -2\tenselem{T}_{1122}^{\R}+4\tenselem{T}_{1212}
\end{bmatrix},\\
\matr{M}^{(1,1)} &=
\begin{bmatrix}
\tenselem{T}_{1111}-\tenselem{T}_{2222} & -\tenselem{T}_{1222}^{\R}-\tenselem{T}_{1112}^\R & -\tenselem{T}_{1222}^{\I}-\tenselem{T}_{1112}^\I\\
-\tenselem{T}_{1222}^{\R}-\tenselem{T}_{1112}^\R & 0 & 0\\
-\tenselem{T}_{1222}^{\I}-\tenselem{T}_{1112}^\I & 0 & 0
\end{bmatrix}.
\end{align*}
\end{example}

It can be seen from \cref{le:sub-problem-general} and the above examples 
that the elementary function $\hij{i}{j}{\matr{U}}$ in \eqref{eq-func-h-new} can always be represented as the sum of a finite number of quadratic forms. 
In fact, for several cases, \emph{e.g.}, the lower order cases \cite[Theorem 4.16]{ULC2019} and the examples in \cref{subsec:glob_conv_examp}, the elementary function $\hij{i}{j}{\matr{U}}$ can be simply represented as a single quadratic form.
In this paper, we also focus on this case, \emph{i.e.}, there exist $1\leq \alpha,\beta\leq\kappa$ such that \begin{equation}\label{eq:quad-form-16}
\hij{i}{j}{\matr{U}}(\matr{\Utwo}) = \vect{z}_{\alpha,\beta}^\T \matr{M}\vect{z}_{\alpha,\beta} + C,
\end{equation}
where $C$ is a constant. 
For the general case in \cref{le:sub-problem-general}, maximizing $\hij{i}{j}{\matr{U}}(\matr{\Utwo})$ may be much more complex, and this will be studied in the future work.

\subsection{Riemannian gradient}
Let $\tilde{g}$ be as in \eqref{definition-f-ext}, and $\ProjGrad{\tilde{g}}{\matr{U}}=\matr{U}\matr{\Lambda}(\matr{U})$ as in \eqref{Rie-gradient-general}.
Then we have the following lemma, which is a direct consequence of \cite[Lemma 3.2]{ULC2019}.
\begin{lemma}\label{lem:ProjGradSubmatrix}
The Riemannian gradient of $\hij{i}{j}{\matr{U}}$ defined in \eqref{eq-func-h-new} at the identity matrix $\matr{I}_2$ is a submatrix of the matrix $\matr{\Lambda}(\matr{U})$:
\begin{align*}\label{eq:gradient-h}
\ProjGrad{\hij{i}{j}{\matr{U}}}{\matr{I}_2}
= \begin{bmatrix}
\matr{\Lambda}(\matr{U})_{ii} & \matr{\Lambda}(\matr{U})_{ij} \\
\matr{\Lambda}(\matr{U})_{ji}  &  \matr{\Lambda}(\matr{U})_{jj}
\end{bmatrix}.
\end{align*}
\end{lemma}

If the elementary function $\hij{i}{j}{\matr{U}}$ has the quadratic form \eqref{eq:quad-form-16}, we further have the following relationship between $\ProjGrad{\hij{i}{j}{\matr{U}}}{\matr{I}_2}$ and the matrix $\matr{M}$ in \eqref{eq:quad-form-16}. For any $1\leq \alpha,\beta\leq\kappa$ instead of $\alpha=2$ and $\beta=1$, \cref{cor:complex gradient-general-1} provides a generalization of \cite[Lemma 4.7]{ULC2019}. Since the proof is similar as in \cite[Lemma 4.7]{ULC2019}, we omit the details here.

\begin{lemma}\label{cor:complex gradient-general-1}
For $1\leq \alpha,\beta\leq\kappa$, we have
\begin{equation*}
\ProjGrad{\hij{i}{j}{\matr{U}}}{\matr{I}_2} =  \alpha 
\begin{bmatrix}
0 & M_{12} + \ui M_{13}\\
-M_{12} + \ui M_{13}& 0
\end{bmatrix}.
\end{equation*}
\end{lemma}

\begin{example}
Let $\tilde{f}$ be as in \eqref{eq:f-tilde-diag}.
Fix $1\leq p\leq d$.
Denote
\begin{equation*}
\tens{V}^{(\ell)}=\tens{A}^{(\ell)} \contr{1} (\matr{U}^{(1)})^{\dagger}\contr{2}\cdots\contr{p-1} (\matr{U}^{(p-1)})^{\dagger}\contr{p+1} (\matr{U}^{(p+1)})^{\dagger}\contr{p+2}\cdots \contr{d} (\matr{U}^{(d)})^{\dagger},
\end{equation*}
where $(\cdot)^{\dagger}$ stands for $(\cdot)^{\H}$ or $(\cdot)^{\T}$. 
Then
\begin{align*}
f_{(p)}(\matr{U})= \sum\limits_{\ell=1}^{L} \alpha_\ell\|\diag{(\tens{W}^{(\ell)})^{(r,r,\cdots,r)}}\|^2  \,\,\,\, \text{with} \,\,\,\,
\tens{W}^{(\ell)}=\tens{V}^{(\ell)}\contr{p} \matr{U}^{\dagger}.
\end{align*}
Let $\vect{v}^{(\ell)}_{p}\eqdef\tens{V}^{(\ell)}_{:,qq\cdots q}\in\CC^{n_p}$.
By \cite[Lemma 5.3]{li2019polar}, we have $\nabla f_{(p)}(\matr{U}) =  2\left[\vect{v}_{1}^{*}\tenselem{W}_{1\cdots 1},\cdots,\vect{v}_{r}^{*}\tenselem{W}_{r\cdots r},\vect{0},\cdots,\vect{0}\right]^\sharp$, 
where $(\cdot)^{\sharp}$ is as in \eqref{eq:def_sharp}. 
It follows from \eqref{Rie-gradient-general} that
{\small
\begin{equation*}\label{eq:Lambda-c-H}
[\matr{\Lambda}^{(p)}(\matr{U})]_{ij} =\left\{
\begin{array}{lcl}
\quad 0,\quad   &  i = j \,\, \mbox{or} \,\, (i,j)\in \mathbb{I}^{(p)}_3,\\
\sum\limits_{\ell=1}^{L} \alpha_\ell \left( \tenselem{W}^{(\ell)}_{j\cdots \p{i}\cdots j}\tenselem{W}_{j\cdots j}^{(\ell)*} - \tenselem{W}_{i\cdots \p{j} \cdots i}^{(\ell)*}\tenselem{W}_{i\cdots i}^{(\ell)} \right)^{\sharp}, & (i,j)\in \mathbb{I}^{(p)}_1,\\
- \sum\limits_{\ell=1}^{L} \alpha_\ell \left( \tenselem{W}^{(\ell)*}_{i\cdots \p{j} \cdots i}\tenselem{W}^{(\ell)}_{i\cdots i} \right)^{\sharp}, & (i,j)\in \mathbb{I}^{(p)}_2, \\
- [\matr{\Lambda}^{(p)}(\matr{U})]^{*}_{ij}, & (j,i)\in \mathbb{I}^{(p)}.
\end{array}
\right.
\end{equation*}}
\end{example}

\section{Jacobi-G algorithm on the unitary group}\label{sec:Jacobi-G algorithm}
\subsection{Jacobi-G algorithm and subproblem} \label{subproblem-general-2}
Let $\tilde{g}$ be an abstract function as in \eqref{definition-f-ext}. The following general \emph{gradient-based Jacobi-type} (Jacobi-G) algorithm was proposed in \cite{ULC2019} to maximize $\tilde{g}$.

\begin{algorithm}
\caption{Jacobi-G algorithm.}\label{Jacobi-G-general}
\begin{algorithmic}[1]
\STATE {\bf Input:} A positive constant $0<\delta\le\sqrt{2}/n$ and an initial point $\matr{U}_{0}$.
\STATE {\bf Output:} Sequence of iterations $\matr{U}_{k}$.
\FOR {$k=1,2,\cdots$ until a stopping criterion is satisfied}
\STATE Choose an index pair $(i_k,j_k)$ satisfying
\begin{equation}\label{inequality-gra-based-2}
\|\ProjGrad{\hij{i_k}{j_k}{\matr{\matr{U}_{k-1}}}}{\matr{I}_2}\| \geq \delta
\|\ProjGrad{\tilde{g}}{\matr{U}_{k-1}}\|.
\end{equation}
\vspace{-0.25cm}
\STATE Find $\matr{\Utwo}_k$ that maximizes $h_k(\matr{\Utwo})\eqdef\hij{i_k}{j_k}{\matr{U}_{k-1}}(\matr{\Utwo})$.
\STATE Update $\matr{U}_k = \matr{U}_{k-1} \Gmat{i_k}{j_k}{\matr{\Utwo}_k}$.
\ENDFOR
\end{algorithmic}
\end{algorithm}

\begin{remark}\label{re:jacobi-G-pair}
According to \cite[Remark 4.2]{ULC2019}, in \cref{Jacobi-G-general}, we can always find an index pair $(i_k,j_k)$ such that the inequality \eqref{inequality-gra-based-2} is satisfied.
\end{remark}
\fin


In this paper, we assume that the function $\tilde{g}$ in \eqref{definition-f-ext} is scale invariant, and $h_k(\matr{\Psi})$ in \cref{Jacobi-G-general} can always be expressed as a quadratic form
\begin{eqnarray}\label{eq:quad-form-12}
h_{k}(\matr{\Psi}) =  \vect{z}_{ \alpha,\beta}^{\T}\matr{M}\vect{z}_{ \alpha,\beta}  + C,
\end{eqnarray}
as in \eqref{eq:quad-form-16} with $1\leq \alpha, \beta \leq \kappa$ and  $\kappa\geq 1$.
Note that $\|  \vect{z}_{ \alpha,\beta}  \|^2=1$.
There may exist many different choices of $\matr{M}$ in \eqref{eq:quad-form-12}. For example, if $\widetilde{\matr{M}} = \matr{M} + C\matr{I}_3$, 
then we also have a quadratic form $h_{k}(\matr{\Psi}) =  \vect{z}_{\alpha,\beta }^{\T}\widetilde{\matr{M}}\vect{z}_{ \alpha,\beta }  $. It is easy to see that $\matr{M}$ and $\widetilde{\matr{M}}$ have the same eigenvector.

\subsection{Convergence analysis}\label{subsec:converg_analy}
 For the symmetric matrix $\matr{M}\in \RR^{3\times 3}$ in \eqref{eq:quad-form-12}, we denote its eigenvalues by $\lambda_{1}\geq\lambda_2\geq\lambda_3$ in a descending order, and their corresponding eigenvectors by $\vect{w}=(w_1,w_2,w_3)^{\T},\vect{u}=(u_1,u_2,u_3)^{\T}$ and $\vect{v}=(v_1,v_2,v_3)^{\T}$.
It is clear that $\max h_k(\matr{\Psi}) =
\lambda_{1} + C$.
Without loss of generality, we always choose $\vect{w}$ such that $w_1\geq 0$. 
Then, we can choose $\theta_k\in [0, \frac{\pi}{2\alpha}]$ such that 
%
\begin{equation}\label{eq:omega_t-2}
\vect{w} =  (\cos \alpha\theta_k,  -\sin \alpha\theta_k\cos\beta\phi_k, -\sin \alpha\theta_k\sin\beta\phi_k)^{\T}\in \mathbb{R}^3
\end{equation} 
by \eqref{eq:omega_t}.
We now present three conditions about the matrix $\matr{M}$ 
as follows:
\begin{enumerate}
\item[(C1)] there exists a universal $\varepsilon>0$ such that $\frac{\lambda_{2}-\lambda_3}{\lambda_1-\lambda_3}\leq 1-\varepsilon$; 
\item[(C2)] $u_1=0$;
\item[(C3)] $v_1=0$.
\end{enumerate}
It can be seen that the above three conditions are all independent of the specific choice of the matrix $\matr{M}$ in \eqref{eq:quad-form-12}. For example, a matrix $\matr{M}$ satisfies the condition (C1) if and only if $\matr{M}+C\matr{I}_3$ satisfies this condition for any constant $C\in\RR$.

 To prove the global convergence of \cref{Jacobi-G-general} under any one of the three conditions (C1), (C2) and (C3), we need a lemma about the relationship between $\|\matr{U}_k - \matr{U}_{k-1}\|$ and $w_1$. 

\begin{lemma}\label{le:global ieq-complex-general-13-3}
Let the matrix $\matr{M}$ be as in \eqref{eq:quad-form-12}, and
$\theta_{k}$ be the corresponding parameter in \eqref{eq:omega_t-2}.
Then 
\begin{align}
\|\matr{U}_k - \matr{U}_{k-1}\| \leq 2\sqrt{1-w_1^2} 
\leq 2\sqrt{2}\sqrt{1-w_1}.\label{eq:thetaalpha-general}
\end{align}
\end{lemma}
\begin{proof}
Note that $\theta_k\in [0, \frac{\pi}{2\alpha}]$. We have
$1-\cos^2 \alpha \theta_k=\sin^2\alpha \theta_k\geq \sin^2\theta_k=1-\cos^2\theta_k.$  
Then, 
%
\begin{align*}
\|\matr{U}_k - \matr{U}_{k-1}\|^2
= \|\Gmat{i_k}{j_k}{\matr{\Psi}_k} - \matr{I}_{n}\|^2= 4(1-\cos\theta_k)
\leq 4(1-\cos^2 \alpha \theta_k)\leq 8(1-w_1),\label{eq:thetaalpha-general-2}
\end{align*}
where the first equation holds since $\matr{U}_k = \matr{U}_{k-1} \Gmat{i_k}{j_k}{\matr{\Utwo}_k}$, and the last inequality follows from $1+w_1 \leq 2$.
The proof is complete.
\end{proof}

\subsubsection{About the condition (C1) for global convergence}
The condition (C1) was in fact first introduced in \cite[Lemma 7.3]{ULC2019}.
As in \cite[Section 7]{ULC2019}, the following result about the global convergence of \cref{Jacobi-G-general} can be easily proved.

\begin{theorem}\label{thm:accumulation_points}
In \cref{Jacobi-G-general}, if the elementary function $h_k(\matr{\Psi})$ in \eqref{inequality-gra-based-2} can always be expressed as a quadratic form in \eqref{eq:quad-form-12}, and there exists a universal $\varepsilon>0$ such that the condition (C1) is always satisfied, then the iterations $\matr{U}_k$ converge to the only limit point $\matr{U}_{*}$. 
\end{theorem}

\begin{proof}
We first prove that
\begin{equation}\label{le:global ieq-complex-general-13-1}
h_{k}(\matr{\Psi}_{k}) - h_{k}(\matr{I}_2) \geq \frac{\varepsilon}{4 \alpha} \|\ProjGrad{h_{k}}{\matr{I}_2}\| \, \|\matr{\Psi}_{k} - \matr{I}_{2}\|
\end{equation}
holds for all $k$. 
In \cite[Equation (7.3)]{ULC2019}, it was shown that
$M_{12}^2 + M_{13}^2 \leq (1-w_1^2)\left((\lambda_1-\lambda_3)^2+(\lambda_2-\lambda_3)^2\right).$ 
By \cref{cor:complex gradient-general-1}, we have
$\|\ProjGrad{h_{k}}{\matr{I}_2}\|^2 = 2 \alpha^2 (M_{12}^2+M_{13}^2)$.  
Together with \cref{le:global ieq-complex-general-13-3}, following the proof of \cite[Lemma 7.3]{ULC2019}, we can check that \eqref{le:global ieq-complex-general-13-1} holds. It is easy to see that $\|\matr{U}_k - \matr{U}_{k-1}\|=\|\matr{\Psi}_{k} - \matr{I}_{2}\|$ from $\matr{U}_k = \matr{U}_{k-1} \Gmat{i_k}{j_k}{\matr{\Utwo}_k}$. 
Then,
based on \cref{theorem-SU15} and \eqref{le:global ieq-complex-general-13-1}, the 
proof is complete.
\end{proof}
\begin{remark}\label{le:global ieq-complex-general-13-1-gamma-2}
It can be seen from \cite[Lemma 7.3]{ULC2019} that, in the quadratic form \eqref{eq:quad-form-12} with $\kappa=2$, if the condition (C1) is always satisfied for a universal positive constant $\varepsilon>0$,
then
\begin{equation}\label{ineq:h_grad-13-1}
h_{k}(\matr{\Psi}_{k}) - h_{k}(\matr{I}_2) \geq \frac{\varepsilon}{8} \|\ProjGrad{h_{k}}{\matr{I}_2}\| \, \|\matr{\Psi}_{k} - \matr{I}_{2}\|
\end{equation}
holds for all $k$.
The equation \eqref{le:global ieq-complex-general-13-1} can be seen as an extension of the above equation \eqref{ineq:h_grad-13-1}.
\end{remark}
\fin

\subsubsection{About the conditions (C2) and (C3) for global convergence}

The conditions (C2) and (C3) are two new ones proposed in this paper.
In this subsection, we mainly prove  the global convergence result of \cref{Jacobi-G-general} under any one of these two conditions. 
 Before that, we first present two examples to explain the condition (C2).

\begin{example}\label{exam:spec-M-A}
If the matrix $\matr{M}$ in the quadratic form \eqref{eq:quad-form-12} has a special structure
\begin{equation}\label{eq:Mmatrix-general}
\matr{M} = \begin{bmatrix}
M_{11} & M_{12} & M_{13} \\
M_{12} & 0 & 0\\
M_{13} & 0 & 0
\end{bmatrix},
\end{equation}
then $\lambda_2=0$ and the second eigenvector $\vect{u}$ satisfies $u_1=0$, and thus the matrix $\matr{M}$ satisfies the condition (C2).
In fact, from the characteristic polynomial of $\matr{M}$, we have 
$\lambda\big(\lambda^2 - M_{11}\lambda - (M_{12}^2+M_{13}^2)\big) = 0.$ 
It follows that
\begin{eqnarray*}\label{lambda-complex-general}
\begin{aligned}
\lambda_{1} = \frac{M_{11}+\sqrt{M_{11}^2+4(M_{12}^2+M_{13}^2)}}{2}
\geq \lambda_2 = 0 \geq \lambda_3 = \frac{M_{11}-\sqrt{M_{11}^2+4(M_{12}^2+M_{13}^2)}}{2}.
\end{aligned}
\end{eqnarray*}
Then, according to $\matr{M}\vect{u}=\lambda_2\vect{u}$, we can easily verify that $u_1=0$. 
 It will be shown in \cref{subsec:glob_conv_examp} that the objective function \eqref{eq:f-tilde-diag} satisfies the structure \eqref{eq:Mmatrix-general}.
\end{example}

\begin{example}
We now present a matrix $\matr{M}$, which satisfies the condition (C2), while does not have the structure \eqref{eq:Mmatrix-general}.
Let
\begin{equation*}\label{eq:Mmatrix-general-2}
\matr{M} = \begin{bmatrix}
1 & -2 &-2 \\
-2 & 1.5 & 1.5 \\
-2 & 1.5 & 1.5
\end{bmatrix}.
\end{equation*}
It can be calculated that $\lambda_1=5,\lambda_2=0$, $\lambda_3=-1$ and
\begin{equation*}\label{eq:Mmatrix-general-3}
[\vect{w}, \vect{u}, \vect{v}] = \begin{bmatrix}
\frac{\sqrt{3}}{3} & 0 &\frac{\sqrt{6}}{3} \\
-\frac{\sqrt{3}}{3} & \frac{\sqrt{2}}{2} & \frac{\sqrt{6}}{6} \\
-\frac{\sqrt{3}}{3} & -\frac{\sqrt{2}}{2} & \frac{\sqrt{6}}{6}
\end{bmatrix}.
\end{equation*}
Therefore, it satisfies (C2), and thus the structure \eqref{eq:Mmatrix-general} is not equivalent to the condition (C2). 
\end{example}
\begin{theorem}\label{thm:accumulation_points-3}
In \cref{Jacobi-G-general}, if the elementary function $h_k(\matr{\Psi})$ in \eqref{inequality-gra-based-2} can always be expressed as a quadratic form in \eqref{eq:quad-form-12}, and the matrix $\matr{M}$ satisfies the condition (C2) or (C3),
then the iterations $\matr{U}_k$ converge to the only limit point $\matr{U}_{*}$.
\end{theorem}

\begin{proof}
We only prove the case for condition (C2). The other case for (C3) is similar. 
Let $\matr{\Psi}_k$ and $\vect{w}$ be the optimal solutions to maximize $h_k(\matr{\Psi})$.
As in equation \eqref{eq:quad-form-12}, we let $\lambda_{1}\geq\lambda_2\geq\lambda_3$ be the eigenvalues of $\matr{M}$ with eigenvectors $\vect{w},\vect{u}$ and $\vect{v}$, respectively.
Then
\begin{equation}\label{eq:lambdaM-general}
\matr{M}=\lambda_{1}\vect{w}\vect{w}^{\T}+\lambda_2\vect{u}\vect{u}^{\T}+\lambda_3\vect{v}\vect{v}^{\T}.
\end{equation}
By the condition (C2), we see that $w_1^2+v_1^2=1$. Then
\begin{align}
\quad \,\, M_{12}^2+M_{13}^2 &= (\lambda_{1}w_1w_2+\lambda_3v_1v_2)^2 + (\lambda_{1}w_1w_3+\lambda_3v_1v_3)^2\notag\\
&= \lambda_{1}^2w_1^{2}(w_2^2+w_3^2) +2\lambda_{1}\lambda_3w_1v_1(w_2v_2+w_3v_3) + \lambda_3^2v_1^2(v_2^2+v_3^2)\notag\\
&= \lambda_{1}^2w_1^{2}(1-w_1^{2}) - 2\lambda_{1}\lambda_3w_1^2v_1^2 + \lambda_3^2v_1^2(1-v_1^2)\notag\\
&= (1-w_1^{2})(\lambda_{1}-\lambda_3)^2(1-v_1^2)\notag\\
&\leq 2(1-w_1)(\lambda_{1}-\lambda_3)^2,\label{ieq:M1213-general}
\end{align}
where the inequality follows from $1+w_1\leq 2$ and $1-v_1^2\leq 1$. Hence, it follows from \cref{cor:complex gradient-general-1} and \eqref{ieq:M1213-general} that
\begin{eqnarray}\label{ieq:gradient lambda13-general}
\begin{aligned}
\|\ProjGrad{h_{k}}{\matr{I}_2}\|^2 = 2 \alpha^2 (M_{12}^2+M_{13}^2)
\leq 4 \alpha^2 (1-w_1)(\lambda_{1}-\lambda_3)^2.
\end{aligned}
\end{eqnarray}
Note that
\begin{align}
h_{k}(\matr{\Psi}_{k}) - h_{k}(\matr{I}_2)
= \vect{w}^{\T}\matr{M}\vect{w} - M_{11}
= \lambda_{1} - (\lambda_{1}w_1^{2}+\lambda_3v_1^2)= (\lambda_{1}-\lambda_3)(1-w_1^{2}),\label{ieq:fh-general}
\end{align}
where the last two equalities result from \eqref{eq:lambdaM-general} and $w_1^{2}+v_1^2= 1$, respectively. From \eqref{eq:thetaalpha-general}, \eqref{ieq:gradient lambda13-general} and \eqref{ieq:fh-general}, we have that 
\begin{align}
h_{k}(\matr{\Psi}_{k}) - h_{k}(\matr{I}_2) &= (\lambda_{1}-\lambda_3)\sqrt{1-w_1}(1+w_1)\sqrt{1-w_1}\notag\\
&\geq \frac{1}{2\alpha} \|\ProjGrad{h_{k}}{\matr{I}_2}\|\sqrt{1-w_1}\geq \frac{\sqrt{2}}{8\alpha} \|\ProjGrad{h_{k}}{\matr{I}_2}\|\|\matr{\Psi}_{k} - \matr{I}_{2}\|.\label{le:global ieq-complex-general}
\end{align}
 Then, by \cref{theorem-SU15}, the iterations $\matr{U}_k$ converge to the only limit point $\matr{U}_{*}$.
The proof is complete.
\end{proof}
\fin

\begin{remark}\label{remark:M-eigen}
For the matrix $\matr{M}$ of form \eqref{eq:Mmatrix-general}, we have the following explicit expression:
\begin{equation*}\label{eq:subproblem solution-complex-general}
\vect{w} =\left\{
\begin{array}{lcl}
\textsf{Normalize}\Big(\textsf{sign}(M_{13})\cdot(\frac{M_{11}+\sqrt{M_{11}^2+4(M_{12}^2+M_{13}^2)}}{2M_{13}}, \frac{M_{12}}{M_{13}},1)^{\T}\Big), & M_{13} \neq 0, \\
\textsf{Normalize}\Big(\textsf{sign}(M_{12})\cdot(\frac{M_{11}+\sqrt{M_{11}^2+4M_{12}^2}}{2M_{12}},1,0)^{\T}\Big), & M_{12} \neq 0 \,\, \mbox{and} \,\, M_{13} = 0,\\
(1,0,0)^{\T} , & M_{12} = 0 \,\, \mbox{and} \,\, M_{13} = 0,
\end{array}
\right.
\end{equation*}
where $\textsf{Normalize}(\bar{\vect{w}})\eqdef\frac{\bar{\vect{w}}}{\|\bar{\vect{w}}\|}$, and the indicator function \textsf{sign}($x$) is defined as
\begin{equation*}
\textsf{sign}(x)=\bigg\{
\begin{aligned}
1,\ \ & x > 0, \\
-1,\ \ & x < 0.
\end{aligned}
\end{equation*}
\end{remark}
\begin{remark}
For real case, the unit vector $\vect{z}_{\alpha,\beta}$ in \eqref{eq:omega_t} will be 
\begin{equation}\label{eq:omega_t-real}
\vect{z}_{\alpha} = (\cos \alpha\theta, -\sin \alpha\theta)^{\T} \in \mathbb{R}^2.
\end{equation} 
In this case, according to the proof of \eqref{le:global ieq-complex-general}, we see that if the elementary function 
can always be expressed as a quadratic form $\vect{z}_{\alpha}^\T \matr{M}\vect{z}_{\alpha} + C$, then the global convergence of real Jacobi-G algorithm on the orthogonal group can be directly obtained without any condition. 
\end{remark}

\begin{remark}\label{remar-com-jac-G}
 We now analyze the computational complexity of \Cref{Jacobi-G-general}. 
Note that the elementary function \eqref{eq:quad-form-12} can be represented as a single quadratic form.
The complexity of solving the subproblem $\max h_k(\matr{\Psi})$ is $\mathcal{O}(1)$, since the solution can be obtained from an eigenvector corresponding to the maximal eigenvalue of $\matr{M} \in \RR^{3\times 3}$. 
At each iteration of the Jacobi-G algorithm, only the elements of the tensor with one of the indices $i$ or $j$ are updated. The computational complexity of per update is $\mathcal{O}(8\kappa n)$. In particular, for problem \eqref{ex:JATD-S}, let $\tens{W}_k$ be the current tensor at $k$-th iteration. If an optimal solution $\matr{\Psi}_k$ has been found, we can modify the current tensor $\tens{W}_k$ by
\begin{equation}\label{eq:contraction-complexity}
\tens{W}_{k+1}= \tens{W}_k\contr{1}\matr{G}(i_k,j_k,\matr{\Psi}_k)^{\dagger} \contr{2}\matr{G}(i_k,j_k,\matr{\Psi}_k)^{\dagger}\cdots\contr{d}\matr{G}(i_k,j_k,\matr{\Psi}_k)^{\dagger}.
\end{equation}
The computational complexity of \eqref{eq:contraction-complexity} is $\mathcal{O}(4d n)$.
Since we need the matrix $\matr{\Lambda}(\matr{U})$ in \eqref{Rie-gradient-general}, the complexity of the gradient inequality is $\mathcal{O}(n^2)$. 
Above all, the computational complexity of per sweep, including $\frac{n(n-1)}{2}$ iterations, is $\mathcal{O}(n^4)$.
\fin
\end{remark}

\section{Jacobi-MG algorithm on the product of unitary groups}\label{sec:Jacobi-MG algorithm}

\subsection{Jacobi-MG algorithm}
In this subsection, we develop a \emph{gradient-based multi-block Jacobi-type} (Jacobi-MG) algorithm on the product of unitary groups $\Upsilon$ to maximize the objective function $\tilde{f}$ in \eqref{definition-f-ext-2}.
As in \cref{subsec:elementary}, for $1\leq p\leq d$,
we define the pair index set $\mathbb{I}^{(p)}\eqdef\{(i,j):1\leq i<j \leq n_p\}$, and divide it into three different subsets
{\small
\begin{equation*}
\mathbb{I}^{(p)}_1=\{(i,j):1\leq i<j \leq r_p\}, \ \mathbb{I}^{(p)}_2=\{(i,j):1\leq i \leq r_p <j \leq n_p\}, \ \mathbb{I}^{(p)}_3=\{(i,j):r_p< i<j \leq n_p\}.
\end{equation*}}
Let $\vect{\upsilon}=(\matr{U}^{(1)},\cdots, \matr{U}^{(d)})\in\Upsilon$ and
\begin{eqnarray*}\label{definition-h-tildeg}
\begin{aligned}
\tilde{g}_{(p)}: \UN{n_p} \longrightarrow \RR, \quad
\matr{U} \longmapsto \tilde{f}(\matr{U}^{(1)}, \cdots, \matr{U}^{(p-1)}, \matr{U}, \matr{U}^{(p+1)}, \cdots, \matr{U}^{(d)})
\end{aligned}
\end{eqnarray*}
be the restricted function defined as in \eqref{definition-h}. 
Let
\begin{align}
h_{(i,j),\vect{\upsilon}}^{(p)}:\UN{2} \longrightarrow \RR,\quad {\matr{\Utwo}}\longmapsto \tilde{g}_{(p)}(\matr{U}^{(p)}\Gmat{i}{j}{\matr{\Utwo}})
\label{eq-func-h-new-2}
\end{align}
be the elementary function of $\tilde{g}_{(p)}$ defined as in \eqref{eq-func-h-new}.
Let $n_{\rm max}\eqdef\max_{1\leq p\leq d}n_p$. 
The general framework of the Jacobi-MG algorithm on $\Upsilon$ can be formulated as in \cref{al:Jacobi-G-general}.

\begin{algorithm}
\caption{Jacobi-MG algorithm.}\label{al:Jacobi-G-general}
\begin{algorithmic}[1]
\STATE {\bf Input:} A positive constant $0 < \delta \leq\sqrt{\frac{2}{dn_{\rm max}(n_{\rm max}-1)}}$ and an initial point $\vect{\upsilon}_{0}$.\\
\STATE {\bf Output:} Sequence of iterations $\vect{\upsilon}_{k}$.
\FOR {$k=1,2,\cdots$ until a stopping criterion is satisfied}
\STATE Choose $p_k$ and an index pair $(i_k,j_k)\in \mathbb{I}^{(p_k)}$ satisfying
\begin{equation}\label{inequality-gra-based}
\|\ProjGrad{h_{(i_k,j_k),\vect{\upsilon}_{k-1}}^{(p_k)}}{\matr{I}_2}\| \geq \delta
\|\ProjGrad{\tilde{f}}{\vect{\upsilon}_{k-1}}\|.
\end{equation}
\vspace{-0.25cm}
\STATE Find $\matr{\Psi}_k$ that maximizes $h_k(\matr{\Psi})\eqdef h_{(i_k,j_k),\vect{\upsilon}_{k-1}}^{(p_k)}(\matr{\Psi})$.
\STATE Update $\matr{U}^{(p_k)}_k = \matr{U}^{(p_k)}_{k-1} \Gmat{i_k}{j_k}{\matr{\Psi}_k}$,
and $\matr{U}^{(p)}_k = \matr{U}^{(p)}_{k-1}$ for $p\neq p_k$.
\STATE Set $\vect{\upsilon}_k = (\matr{U}_{k}^{(1)},\cdots,\matr{U}_{k}^{(d)})$.
\ENDFOR
\end{algorithmic}
\end{algorithm}

 As the inequality \eqref{inequality-gra-based-2} in \cref{Jacobi-G-general}, the above inequality \eqref{inequality-gra-based} restricted to the $p_k$-th block is to guarantee the global convergence result \Cref{th:global convergence-general}. 
Now we first prove that \Cref{al:Jacobi-G-general} is well defined.
In other words, we can always find a pair $(i_k,j_k)$ and mode $p_k$ such that the inequality \eqref{inequality-gra-based} is satisfied. 
\begin{lemma}\label{le:hf-complex-general}
Let $\tilde{f}$ and $h^{(p)}_{(i,j),\vect{\upsilon}}$ be defined as in \eqref{definition-f-ext-2} and \eqref{eq-func-h-new-2}, respectively. Then there always exist $1\leq p\leq d$ and $(i,j)\in \mathbb{I}^{(p)}$ such that
\begin{equation*}\label{ieq:al}
\|\ProjGrad{h^{(p)}_{(i,j),\vect{\upsilon}}}{\matr{I}_2}\| \geq \delta \|\ProjGrad{\tilde{f}}{\vect{\vect{\upsilon}}}\|,
\end{equation*}
where $0 < \delta\leq\sqrt{\frac{2}{dn_{\rm max}(n_{\rm max}-1)}}$ is a fixed positive constant.
\end{lemma}
\begin{proof}
Let $\ProjGrad{\tilde{g}_{(p)}}{\matr{U}}=\matr{U}\matr{\Lambda}^{(p)}(\matr{U})$ as in \eqref{Rie-gradient-general}.
By \cref{lem:ProjGradSubmatrix}, we have
\begin{align*}
\sum_{(i,j)\in \mathbb{I}^{(p)}}\|\ProjGrad{h^{(p)}_{(i,j),\vect{\upsilon}}}{\matr{I}_2}\|^2\geq \|\matr{\Lambda}^{(p)}(\matr{U}^{(p)})\|^2= \|\ProjGrad{\tilde{g}_{(p)}}{\matr{U}^{(p)}}\|^2.
\end{align*}
It follows that
\begin{align*}
\frac{d}{2}n_{\rm max}(n_{\rm max}-1)\max_{(i,j)\in\mathbb{I}^{(p)},1\leq p\leq d}\|\ProjGrad{h^{(p)}_{(i,j),\vect{\upsilon}}}{\matr{I}_2}\|^2&\geq \sum_{p=1}^{d}\sum_{(i,j)\in \mathbb{I}^{(p)}}\|\ProjGrad{h^{(p)}_{(i,j),\vect{\upsilon}}}{\matr{I}_2}\|^2\\
&\geq \sum_{p=1}^{d}\|\ProjGrad{\tilde{g}_{(p)}}{\matr{U}^{(p)}}\|^2= \|\ProjGrad{\tilde{f}}{\vect{\vect{\upsilon}}}\|^2.
\end{align*}
The proof is complete. 
\end{proof}

\begin{remark}\label{re:Jacobi-C}
In \Cref{al:Jacobi-G-general}, a more natural way of choosing the block $p_k$ and index pair $(i_k,j_k)$ is according to a cyclic ordering. For example, we first choose $1\leq p_k\leq d$ cyclically, and then for each $p_k$, we choose the index pair in a cyclic way as follows:
\begin{equation*}\label{partial-cyclic-1}
\begin{split}
(1,2) \to (1,3) \to \cdots \to (1,n_k) \to 
 (2,3) \to \cdots \to (2,n_k) \to 
 \cdots \to (n_k-1,n_k) \to 
(1,2) \to (1,3) \to \cdots.
\end{split}
\end{equation*}
In this way, we call it a \emph{Jacobi-MC algorithm} on the product of unitary groups. 
This cyclic choice of index pair $(i_k,j_k)$ was used in \cite{martin2008jacobi}.
\end{remark}

\subsection{Global convergence}
As in \cref{subproblem-general-2}, we now assume that the function $\tilde{f}$ in \eqref{definition-f-ext-2} is scale invariant, and $h_k(\matr{\Psi})$ in \cref{al:Jacobi-G-general} can always be expressed as a quadratic form in \eqref{eq:quad-form-12}.
 By a similar method for proving \eqref{le:global ieq-complex-general-13-1} and \eqref{le:global ieq-complex-general}, we can get the following result.

\begin{lemma}\label{le:global ieq-complex-general-13-12}
If the quadratic form \eqref{eq:quad-form-12} satisfies the condition (C1) or (C2) or (C3),
then 
\begin{equation*}\label{ineq:h_grad-13-12}
\tilde{f}(\vect{\upsilon}_k)-\tilde{f}(\vect{\upsilon}_{k-1}) \geq \eta \|\ProjGrad{\tilde{f}}{\vect{\upsilon}_{k-1}}\| \|\vect{\upsilon}_{k} - \vect{\upsilon}_{k-1}\|
\end{equation*}
holds for all $k$, where $\eta=\frac{\varepsilon \delta}{4 \alpha}$ under the condition (C1) and $\eta=\frac{\sqrt{2}\delta}{8 \alpha }$ under the condition (C2) or (C3). 
\end{lemma}

Now, we apply \cref{theorem-SU15} and \cref{le:global ieq-complex-general-13-12} to get the global convergence of \cref{al:Jacobi-G-general} easily.
\begin{theorem}\label{th:global convergence-general}
In \cref{al:Jacobi-G-general} for the objective function $\tilde{f}$ in \eqref{definition-f-ext-2}, if the quadratic form \eqref{eq:quad-form-12} always satisfies any one of the conditions (C1), (C2) and (C3), then for any initial point $\vect{\upsilon}_0$, the iterations $\vect{\upsilon}_k$ converge to the only limit point $\vect{\upsilon}_{*}$.
\end{theorem}

\subsection{Joint approximate tensor diagonalization}\label{subsec:glob_conv_examp}

In this subsection, we show that the condition (C2) is always satisfied when $\tilde{f}$ is the function \eqref{eq:f-tilde-diag}. 
Then, by \cref{th:global convergence-general}, we get the global convergence of \cref{al:Jacobi-G-general} for this function. 

\begin{lemma}\label{lem:M-JATD-tilde-2}
Let $\tilde{f}$ be as in \eqref{eq:f-tilde-diag}, and $h_{(i,j),\vect{\upsilon}}^{(p)}$ be as in \eqref{eq-func-h-new-2}.
Then $h_{(i,j),\vect{\upsilon}}^{(p)}$ can always be expressed as a quadratic form \eqref{eq:quad-form-12} with $\alpha=\beta=1$.
Let
\begin{equation}\label{eq:varrho}
\varrho=\bigg\{
\begin{aligned}
1,\ \ & \textrm{if}\ (\cdot)^{\dagger}=(\cdot)^{\H}, \\
-1,\ \ & \textrm{if}\ (\cdot)^{\dagger}=(\cdot)^{\T}.
\end{aligned}
\end{equation}
(i) If $(i,j)\in\mathbb{I}^{(p)}_1$, then $\matr{M}$ has the form \eqref{eq:Mmatrix-general} with
\begin{align*}
M_{11} &= \sum\limits_{\ell=1}^{L} \alpha_\ell \left( |\tenselem{W}_{i\cdots i}^{(\ell)}|^2 + |\tenselem{W}_{j \cdots j}^{(\ell)}|^2 - |\tenselem{W}_{i\cdots \p{j}\cdots i}^{(\ell)}|^2 - |\tenselem{W}_{j \cdots \p{i}\cdots j}^{(\ell)}|^2 \right),\\
M_{12} &= \sum\limits_{\ell=1}^{L} \alpha_\ell \left(\tenselem{W}_{j \cdots j}^{(\ell)\R}\tenselem{W}_{j \cdots \p{i}\cdots j}^{(\ell)\R} + \tenselem{W}_{j \cdots j}^{(\ell)\I}\tenselem{W}_{j \cdots \p{i}\cdots j}^{(\ell)\I}-\tenselem{W}_{i\cdots i}^{(\ell)\R}\tenselem{W}_{i\cdots \p{j}\cdots i}^{(\ell)\R} - \tenselem{W}_{i\cdots i}^{(\ell)\I}\tenselem{W}_{i\cdots \p{j}\cdots i}^{(\ell)\I} \right),\\
M_{13} &= \sum\limits_{\ell=1}^{L} \varrho\alpha_\ell \left(\tenselem{W}_{i\cdots i}^{(\ell)\R}\tenselem{W}_{i\cdots \p{j}\cdots i}^{(\ell)\I} +\tenselem{W}_{j \cdots j}^{(\ell)\R}\tenselem{W}_{j \cdots \p{i}\cdots j}^{(\ell)\I}-\tenselem{W}_{i\cdots i}^{(\ell)\I}\tenselem{W}_{i\cdots \p{j}\cdots i}^{(\ell)\R} - \tenselem{W}_{j \cdots j}^{(\ell)\I}\tenselem{W}_{j \cdots \p{i}\cdots j}^{(\ell)\R} \right).
\end{align*}%
(ii) If $(i,j)\in \mathbb{I}^{(p)}_2$, then $\matr{M}$ has the form \eqref{eq:Mmatrix-general} with
\begin{eqnarray*}
\begin{aligned}
M_{11} &= \sum\limits_{\ell=1}^{L} \alpha_\ell \left(|\tenselem{W}_{i\cdots i}^{(\ell)}|^2 - |\tenselem{W}_{i\cdots \p{j}\cdots i}^{(\ell)}|^2\right),\ 
M_{12} = \sum\limits_{\ell=1}^{L} \alpha_\ell \left( -\tenselem{W}_{i\cdots i}^{(\ell)\R}\tenselem{W}_{i\cdots \p{j}\cdots i}^{(\ell)\R} - \tenselem{W}_{i\cdots i}^{(\ell)\I}\tenselem{W}_{i\cdots \p{j}\cdots i}^{(\ell)\I} \right),\\
M_{13} &= \sum\limits_{\ell=1}^{L} \varrho\alpha_\ell \left( \tenselem{W}_{i\cdots i}^{(\ell)\R}\tenselem{W}_{i\cdots \p{j}\cdots i}^{(\ell)\I} - \tenselem{W}_{i\cdots i}^{(\ell)\I}\tenselem{W}_{i\cdots \p{j}\cdots i}^{(\ell)\R} \right).
\end{aligned}
\end{eqnarray*}
\end{lemma}
By \cref{exam:spec-M-A}, \cref{th:global convergence-general} and \cref{lem:M-JATD-tilde-2}, we can directly get the following result about \cref{al:Jacobi-G-general} for the objective function $\tilde{f}$ in \eqref{eq:f-tilde-diag}.
\begin{corollary}\label{th:global convergence-general-123}
In \cref{al:Jacobi-G-general} for the objective function $\tilde{f}$ in \eqref{eq:f-tilde-diag}, for any initial point $\vect{\upsilon}_0$, the iterations $\vect{\upsilon}_k$ converge to the only limit point $\vect{\upsilon}_{*}$.
\end{corollary}
\begin{remark}\label{re:martin08-remark}
(i) A Jacobi-type algorithm was proposed in \cite{martin2008jacobi} to solve the real order-$3$ tensor case of JATD problem, in which the subproblem is to maximize the sum of squares of the diagonal elements, or the trace.
In this algorithm, each iteration includes three sweeps, and each sweep corresponds to one orientation (the order-$3$ tensor has three different orientations as a cube).\\ 
(ii) It was mentioned in \cite[Section 11]{martin2008jacobi} that ``\emph{the question of the general order-$d$ tensor case is still being explored}". 
 For any order-$d$ tensors, the \cref{al:Jacobi-G-general} chooses one orientation at each iteration, and the subproblem for JATD problem can always be represented as a quadratic form using the unit vector in \eqref{eq:omega_t-real}, in which case the global solution can be calculated. 
Therefore, in some sense, this issue in \cite[Section 11]{martin2008jacobi} has been addressed, and the complex tensor case is considered as well. \\
(iii) Moreover, analogous with the discussions above, it is easy to check that our theoretical results in \cref{th:global convergence-general-123} also directly apply to the \emph{real gradient-based multi-block Jacobi-type algorithm} on the product of orthogonal groups for real JATD problem. We omit the details here. 
\end{remark}

\begin{remark}\label{remar-com-jac-MG}
As in \Cref{remar-com-jac-G}, the computational complexity of per sweep for \Cref{al:Jacobi-G-general} is $\mathcal{O}(dn_{\max}^4)$. In particular, for problem \eqref{ex:JATD}, $\matr{M}$ has a special structure \eqref{eq:Mmatrix-general}, and thus we always have an explicit expression for the solution of the subproblem by \cref{remark:M-eigen}. 
Compared with \Cref{Jacobi-G-general}, the difference is that, since $
\tens{W}_{k+1}= \tens{W}_k\contr{p_k}\matr{G}(i_k,j_k,\matr{\Psi}_k)^{\H}$ at each iteration,
the computational complexity of per update is $\mathcal{O}(4n_{p_{k}})$. In \Cref{al:Jacobi-G-general}, although every $p_k$ and index pair $(i_k,j_k)$ satisfying the inequality \eqref{inequality-gra-based} will lead to a subproblem, it has an analytical solution, and we can compute it very fast. Therefore, \Cref{al:Jacobi-G-general} is of great potential in general. 
\end{remark}
\fin

\section{Jacobi-GP and Jacobi-MGP algorithms}\label{sec:Jacobi-GP algorithm}

\subsection{Jacobi-GP algorithm}

Let $\epsilon>0$ be a small positive constant.
Let $\vect{e}_1=(1,0,0)^{\T}\in\RR^{3}$ be the standard unit vector. 
For the objective function \eqref{definition-f-ext},
we define a new function
\begin{align}\label{eq:h-tilde-M}
\hijtilde{i}{j}{\matr{U}}(\matr{\Utwo}) \eqdef \hij{i}{j}{\matr{U}}(\matr{\Utwo}) - \epsilon\| \vect{z}_{ \alpha, \beta} -\vect{e}_1\|^2,
\end{align}
where $\vect{z}_{ \alpha, \beta}$ is related to $\matr{\Utwo}$ by \eqref{unitary-para-2} and \eqref{eq:omega_t}. 
Now, in \cref{Jacobi-GV-general}, we propose a proximal variant of the Jacobi-G algorithm, which is called the \emph{Jacobi-GP} algorithm.
\begin{algorithm}
\caption{Jacobi-GP algorithm.}\label{Jacobi-GV-general}
\begin{algorithmic}[1]
\STATE {\bf Input:} Two positive constants $0<\delta\le\sqrt{2}/n$ and $\epsilon>0$, an initial point $\matr{U}_{0}$.
\STATE {\bf Output:} Sequence of iterations $\matr{U}_{k}$.
\FOR{ $k=1,2,\cdots$ until a stopping criterion is satisfied}
\STATE Choose an index pair $(i_k,j_k)$ such that the inequality \eqref{inequality-gra-based-2} is satisfied.
\STATE Find $\matr{\Utwo}_k$ that maximizes $\tilde{h}_k(\matr{\Utwo})\eqdef\hijtilde{i_k}{j_k}{\matr{U}_{k-1}}(\matr{\Utwo})$ defined in \eqref{eq:h-tilde-M}.
\STATE Update $\matr{U}_k = \matr{U}_{k-1} \Gmat{i_k}{j_k}{\matr{\Utwo}_k}$.
\ENDFOR
\end{algorithmic}
\end{algorithm}
\begin{remark}\label{re:jacobi-gp}
It is worth noting that inequality \eqref{inequality-gra-based-2} should be satisfied for the original $\hij{i_k}{j_k}{\matr{U}_{k-1}}(\matr{I}_2)$ in step 4 of \cref{Jacobi-GV-general}, and thus, from \cref{re:jacobi-G-pair}, we can always find an index pair $(i_k,j_k)$ such that it is satisfied.
\end{remark}
\fin

\begin{lemma}\label{le:global ieq-complex-general-13}
Let $\hijtilde{i}{j}{\matr{U}}$ be as in \eqref{eq:h-tilde-M}, and $\matr{\Utwo}_{*}$ be the maximizer. Let $(\vect{z}_{\alpha, \beta})_{*}$ be the unit vector $\vect{z}_{\alpha, \beta}$ corresponding to $\matr{\Utwo}_{*}$.
 Then
\begin{align*}
\hij{i}{j}{\matr{U}}(\matr{\Utwo}_{*})-\hij{i}{j}{\matr{U}}(\matr{I}_2)\geq \epsilon\|(\vect{z}_{\alpha, \beta})_{*}-\vect{e}_1\|^2.
\end{align*}
\end{lemma}

\begin{proof}
If $\matr{\Utwo}_{*}$ is the maximizer of $\hijtilde{i}{j}{\matr{U}}$ in \eqref{eq:h-tilde-M}, we see that
\begin{align*}
\hij{i}{j}{\matr{U}}(\matr{\Utwo}_{*}) - \hij{i}{j}{\matr{U}}(\matr{I}_2) &= \hijtilde{i}{j}{\matr{U}}(\matr{\Utwo}_{*}) - \hijtilde{i}{j}{\matr{U}}(\matr{I}_2) +\epsilon\|( \vect{z}_{ \alpha, \beta } )_{*}-\vect{e}_1\|^2\geq \epsilon\|( \vect{z}_{ \alpha, \beta } )_{*}-\vect{e}_1\|^2,
\end{align*}
and thus the proof is complete.
\end{proof}

\begin{corollary}\label{cor:global ieq-complex-general-13}
In \cref{Jacobi-GV-general} for the objective function \eqref{definition-f-ext}, for any initial point $\matr{U}_0$, the iterations $\matr{U}_k$ converge to a stationary point $\matr{U}_{*}$. 
\end{corollary}

\begin{proof}
By \cref{theorem-SU15}, it is sufficient to prove that
\begin{equation}\label{eq:thm2.3-glob}
h_{k}(\matr{\Psi}_{k}) - h_{k}(\matr{I}_2) \geq \epsilon' \|\ProjGrad{h_{k}}{\matr{I}_2}\| \, \|\matr{\Psi}_{k} - \matr{I}_{2}\|
\end{equation}
for a universal constant $\epsilon'>0$ and any $k\geq 1$.
By \cite[Lemma 7.2]{ULC2019}, there exists a universal constant $\tau >0$ such that
\begin{equation}\label{ieq:lemma 7.2}
\tau \|\ProjGrad{h_{k}}{\matr{I}_2}\| \leq \|\matr{\Psi}_{k} - \matr{I}_{2}\|.
\end{equation}
It follows from \eqref{eq:thetaalpha-general} that
\begin{equation}\label{ieq:zgamma-U}
\|( \vect{z}_{\alpha, \beta} )_{*}-\vect{e}_1\|^2 = 2(1-( \vect{z}_{\alpha, \beta} )_{*1}) \geq \frac{1}{4}\|\matr{\Psi}_{k} - \matr{I}_{2}\|^2.
\end{equation}
By \cref{le:global ieq-complex-general-13}, and using \eqref{ieq:lemma 7.2} and \eqref{ieq:zgamma-U}, we have that
\begin{equation}\label{ieq:Jacobi-GP-KL}
h_{k}(\matr{\Psi}_{k}) - h_{k}(\matr{I}_2)  \geq \epsilon\|( \vect{z}_{\alpha, \beta} )_{*}-\vect{e}_1\|^2
\geq \frac{\epsilon}{4}\|\matr{\Psi}_{k} - \matr{I}_{2}\|^2 
\geq \frac{\epsilon\tau}{4}\|\ProjGrad{h_{k}}{\matr{I}_2}\| \, \|\matr{\Psi}_{k} - \matr{I}_{2}\|,
\end{equation}
and then the inequality \eqref{eq:thm2.3-glob} holds by denoting $\epsilon'=\frac{\epsilon\tau}{4}$. 
Moreover, we see that $\matr{U}_*$ is a stationary point from \eqref{ieq:lemma 7.2} and \eqref{ieq:Jacobi-GP-KL}. The proof is complete. 
\end{proof}

\begin{remark}
In some sense, \cref{Jacobi-GV-general} can be seen as a generalization of the Jacobi-PC algorithm proposed in \cite{LUC2017globally}, and the JLROA-GP algorithm proposed in \cite{li2019jacobi}.   
\end{remark}

\subsection{Solving the subproblem}

For the objective function  \eqref{definition-f-new-32}, we assume that the elementary function $\hij{i}{j}{\matr{U}}$ in \eqref{eq-func-h-new} can always be expressed as a quadratic form in \eqref{eq:quad-form-16}.
In this case, we have
\begin{align}\label{eq:h-tilde-subproblem}
\hijtilde{i}{j}{\matr{U}}(\matr{\Utwo}) =  \vect{z}_{ \alpha, \beta } ^\T \matr{M} \vect{z}_{\alpha, \beta}  - \epsilon\| \vect{z}_{\alpha, \beta } -\vect{e}_1\|^2 + C
=  \vect{z}_{\alpha, \beta} ^{\T} \matr{M} \vect{z}_{\alpha, \beta}  + 2\epsilon\langle \vect{z}_{\alpha, \beta} , \vect{e}_1\rangle -2\epsilon+C. 
\end{align}
Then the subproblem \eqref{eq:h-tilde-subproblem} can be rewritten as
\begin{equation}\label{contrained-least-equares}
\vect{z}_{\alpha, \beta}=\argmin_{\|\vect{z}_{\alpha, \beta}\|=1} \left\{-\vect{z}_{\alpha, \beta}^{\T}\matr{M}\vect{z}_{\alpha, \beta} - 2\epsilon\langle \vect{z}_{\alpha, \beta} , \vect{e}_1\rangle +2\epsilon-C \right\},
\end{equation}
which is a classical \emph{constrained least squares} problem. The
problem \eqref{contrained-least-equares} can be solved using the approach in \cite{Gander1988A}, where the authors \cite{Gander1988A} investigated the solvability of the constrained least squares problem based on quadratic eigenvalue problem. 

\subsection{Jacobi-MGP algorithm}

Let $\epsilon>0$ be a small positive constant.
For the objective function \eqref{definition-f-ext-2},
we similarly define a new function
\begin{align}\label{eq:h-tilde-M-P1-gg}
\tilde{h}_{(i,j),\vect{\upsilon}}^{(p)}(\matr{\Utwo}) \eqdef h_{(i,j),\vect{\upsilon}}^{(p)}(\matr{\Utwo}) - \epsilon\| \vect{z}_{\alpha, \beta} -\vect{e}_1\|^2.
\end{align}
Now, in \cref{Jacobi-GV-general-gg}, we propose a proximal variant of the Jacobi-MG algorithm, which is called the \emph{Jacobi-MGP} algorithm.

\begin{algorithm}
\caption{Jacobi-MGP algorithm.}\label{Jacobi-GV-general-gg}
\begin{algorithmic}[1]
\STATE {\bf Input:} Two positive constants $0 < \delta \leq\sqrt{\frac{2}{dn_{\rm max}(n_{\rm max}-1)}}$ and $\epsilon$, an initial point $\vect{\upsilon}_{0}$.\\
\STATE {\bf Output:} Sequence of iterations $\vect{\upsilon}_{k}$.
\FOR {$k=1,2,\cdots$ until a stopping criterion is satisfied}
\STATE Choose $p_k$ and an index pair $(i_k,j_k)\in \mathbb{I}^{(p_k)}$ satisfying
the inequality \eqref{inequality-gra-based} is satisfied.
\STATE Find $\matr{\Psi}_k$ that maximizes $\tilde{h}_k(\matr{\Psi})\eqdef \tilde{h}_{(i_k,j_k),\vect{\upsilon}_{k-1}}^{(p_k)}(\matr{\Psi})$ in \eqref{eq:h-tilde-M-P1-gg}.
\STATE Update $\matr{U}^{(p_k)}_k = \matr{U}^{(p_k)}_{k-1} \Gmat{i_k}{j_k}{\matr{\Psi}_k}$,
and $\matr{U}^{(p)}_k = \matr{U}^{(p)}_{k-1}$ for $p\neq p_k$.
\STATE Set $\vect{\upsilon}_k = (\matr{U}_{k}^{(1)},\cdots,\matr{U}_{k}^{(d)})$.
\ENDFOR
\end{algorithmic}
\end{algorithm}

Analogous with \cref{le:global ieq-complex-general-13} and \cref{cor:global ieq-complex-general-13}, from \cref{theorem-SU15}, we can easily get the following convergence result.
\begin{corollary}
In \cref{Jacobi-GV-general-gg} for the objective function \eqref{definition-f-ext-2}, for any initial point $\vect{\upsilon}_0$, the iterations $\vect{\upsilon}_k$ converge to a stationary point $\vect{\upsilon}_{*}$.
\end{corollary}

\section{Numerical experiments}\label{sec:experiments}

In this section, we report some numerical results using several random examples. 
All the computations are done using MATLAB R2019a and the Tensor Toolbox version 2.6 \cite{TTB2015Software}.
The numerical experiments are conducted on a PC with an Intel\textsuperscript{\textregistered} Core\textsuperscript{TM} i5 CPU at 2.11 GHz and 8.00 GB of RAM in 64bt Windows operation system.
The Jacobi-type algorithms\footnote{ The MATLAB codes are available at \url{https://github.com/mashengz/Jacobi-type}.} we run in this section include the \emph{Jacobi-MG, Jacobi-MC (cyclic), Jacobi-G, Jacobi-GP} and \emph{Jacobi-C}, while the first-order Riemannian optimization methods, including the \emph{steepest descent method, conjugate gradient method} and \emph{Barzilai-Borwein method}, are implemented in the \textsf{manopt}\footnote{This package was downloaded from \url{https://www.manopt.org}.} package \cite{JMLR:v15:boumal14a} as comparisons. 
We choose $\matr{U}_0 =\matr{I}_{n}$ or $\vect{\omega}_0 = (\matr{I}_{n_1},\cdots,\matr{I}_{n_d})\in\Omega$ as the initial point.
In the Jacobi-type algorithms, the stopping criterion is that the sweep number is greater than the maximum sweep number $\emph{maxsweep}=200$.  
The first-order Riemannian optimization methods in the \textsf{manopt} package use the default values of parameters
except the stopping criterion was changed. 
These first-order Riemannian optimization algorithms run up to $\emph{maxiter}=1,000$ iterations.

Let $\tens{A}\in\CC^{n_1\times n_2\times\cdots\times n_d}$.
 We define the ratio 
$\texttt{Rat}\eqdef \sum_{i=1}^{n_{\rm min}}\frac{|\tenselem{A}_{ii\cdots i}|^2}{\|\tens{A}\|^2}$
as a measure for the approximate tensor diagonalization.
Obviously, \texttt{Rat} can be used to measure how much of the elements are condensed on the diagonal of $\tens{A}$. It is equal to 1 if a tensor is diagonal, and in general it is between 0 and 1. If an approximate tensor diagonalization is better, the ratio \texttt{Rat} will be closer to 1.

\begin{figure}[htb!]
\begin{minipage}[b]{0.49\linewidth}
\centering
\includegraphics[height=5.5cm]{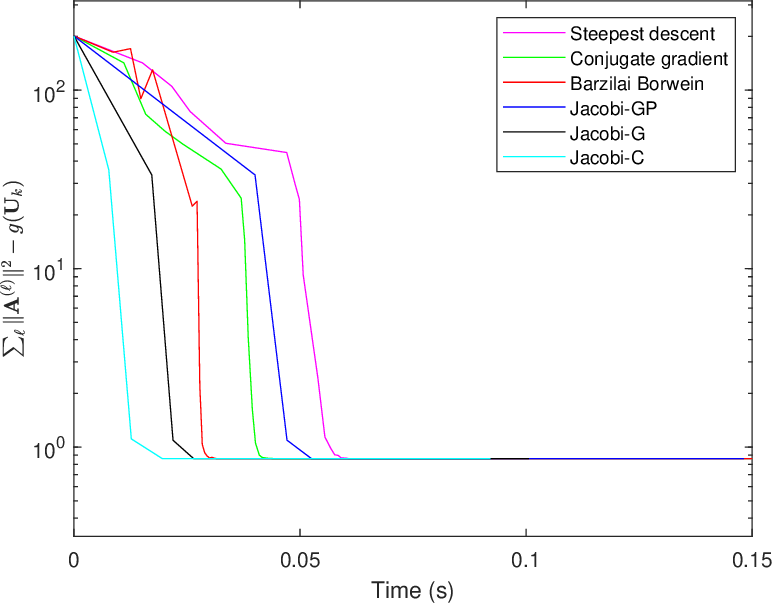}
\end{minipage}
\hfill
\begin{minipage}[b]{0.49\linewidth}
\centering
\includegraphics[height=5.5cm]{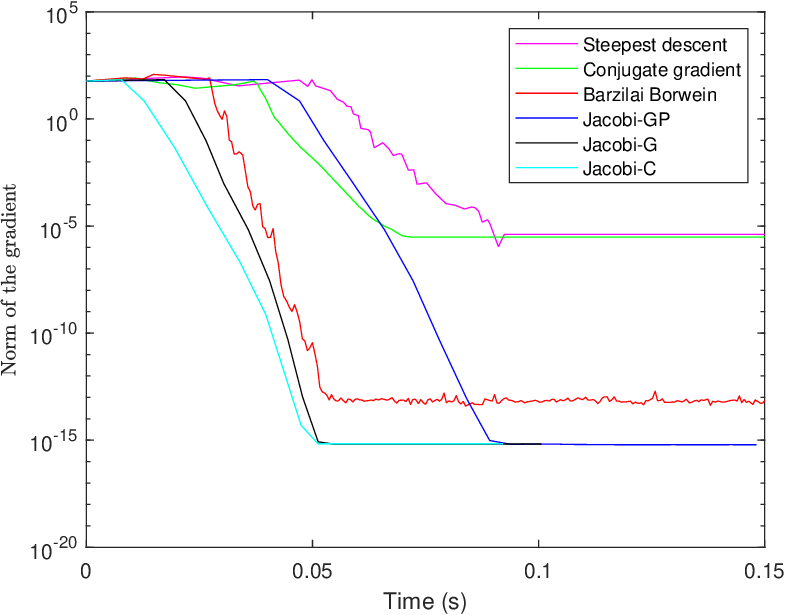}
\end{minipage}
\caption{Experimental results for \cref{ex:3}} \label{fig:ex-joint-diag-fun-grad}
\end{figure}

\begin{example}\label{ex:3}
We consider the JATD-S problem in \eqref{cost-fn-general-1-s} with $n=10$, $r=10$ and $L=8$. Let $d=2$ and $\alpha_\ell=1$ for $1\leq \ell \leq L$. 
Let the matrices $\matr{A}^{(\ell)} = \matr{D}^{(\ell)} \contr{1} \matr{U} \contr{2} \matr{U} + \frac{\matr{E}^{(\ell)}}{\|\matr{E}^{(\ell)}\|} $ for $1\leq \ell \leq L$, where $\tens{D}^{(\ell)}=\textsf{diag}(2*\textsf{randn}(n,1))$ are randomly generated for $1\leq \ell \leq L$, $\matr{U} \in \CC^{n\times n}$ is a random unitary matrix and $\matr{E}^{(\ell)}=2*\textsf{randn}(n,n)$ for $1\leq \ell \leq L$. Here, we have used the notations $\textsf{diag}$ and $\textsf{randn}$, which are the MATLAB built-in functions. It can be calculated that the initial ratio $\texttt{Rat} = 0.3146$.
We choose the parameters $\delta = \frac{\sqrt{2}}{10n}$ and
$\epsilon=10^{-3}$ in Jacobi-G and Jacobi-GP. 
In \cref{fig:ex-joint-diag-fun-grad}, we plot the results of  $\sum_{\ell}\|\matr{A}^{(\ell)}\|^2-g(\matr{U}_k)$ and $\|\ProjGrad{g}{\matr{U}_k}\|$ with respect to the running time. 
In particular, all the algorithms in \cref{fig:ex-joint-diag-fun-grad} arrive at the same ratio $\texttt{Rat}=0.9971$.
\end{example}

\begin{figure}[htb!]
\begin{minipage}[b]{0.49\linewidth}
\centering
\includegraphics[height=5.5cm]{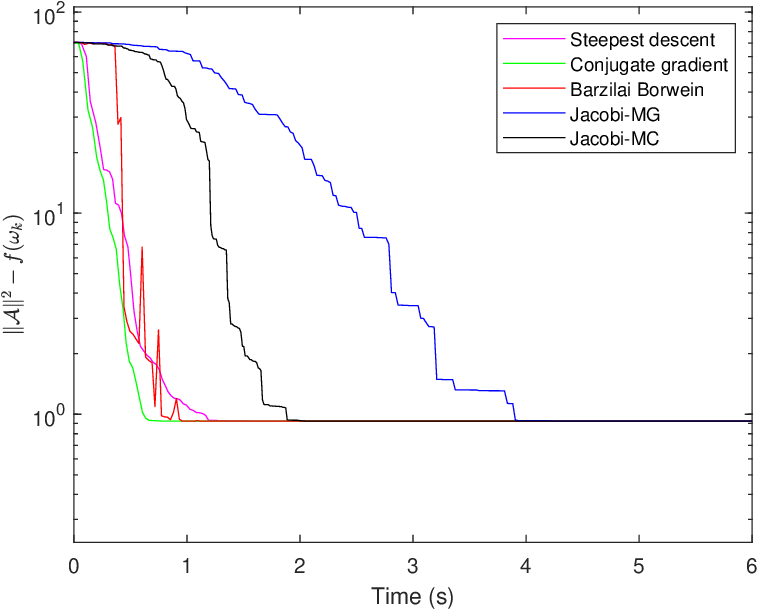}
\end{minipage}
\hfill
\begin{minipage}[b]{0.49\linewidth}
\centering
\includegraphics[height=5.5cm]{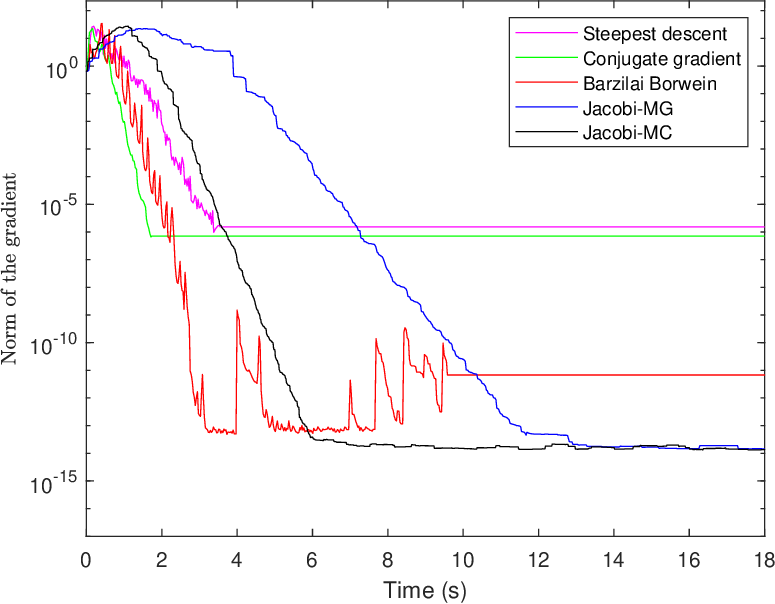}
\end{minipage}
\caption{Experimental results for \cref{ex:2}} \label{fig:ex-complex-diag-fun-grad}
\end{figure}

\begin{example}\label{ex:2}
We consider the JATD problem in \eqref{cost-fn-general-1} with $d=4$, $n_1=n_2=5$, $n_3=n_4=6$ and $r=5$. 
Let $L=1$ and $\alpha_1=1$. 
Let $\tens{A} = \tens{D} \contr{1} \matr{U}^{(1)} \contr{2} \matr{U}^{(2)} \contr{3} \matr{U}^{(3)} \contr{4} \matr{U}^{(4)} + \frac{\tens{E}}{\|\tens{E}\|}$, where $\tens{D}$ is a diagonal tensor with $\tenselem{D}_{jjjj}=\sqrt{j}+j\ui$ for $1\leq j\leq n_{\rm min}=5$, $\matr{U}^{(p)} \in \CC^{n_p\times n_p}$, $p=1,2,3,4$, are random unitary matrices, and $\tens{E}=\textsf{randn}(n_1,n_2,n_3,n_4)+\ui*\textsf{randn}(n_1,n_2,n_3,n_4)$. 
It can be calculated that the initial ratio $\texttt{Rat} = 0.0026$. 
In Jacobi-MG, we choose the parameter $\delta=\frac{1}{10}\sqrt{\frac{2}{d n_{\rm max}(n_{\rm max}-1)}}$. 
The results of  $\|\tens{A}\|^2-f( \vect{\omega}_k)$ and $\|\ProjGrad{f}{\vect{\omega}_k}\|$ with respect to the running time are shown in \Cref{fig:ex-complex-diag-fun-grad}. 
In particular, all the algorithms in \Cref{fig:ex-complex-diag-fun-grad} arrive at the same ratio $\texttt{Rat}=0.9870$.
\end{example}

\begin{figure}[htb!]
\begin{minipage}[b]{0.49\linewidth}
\centering
\includegraphics[height=5.5cm]{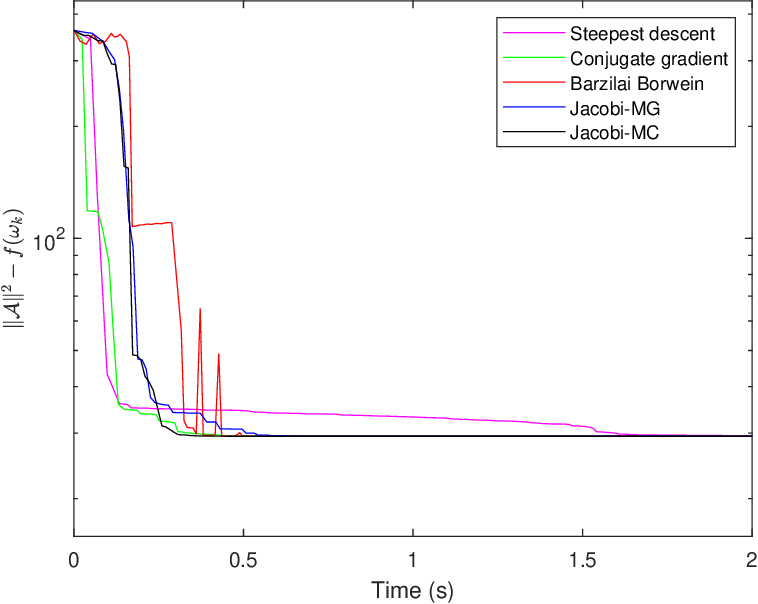}
\end{minipage}
\hfill
\begin{minipage}[b]{0.49\linewidth}
\centering
\includegraphics[height=5.5cm]{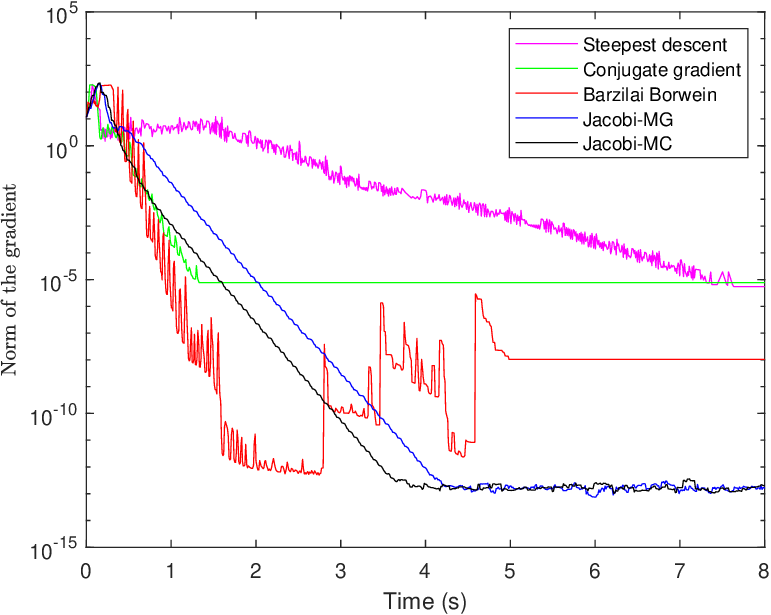}
\end{minipage}
\caption{Experimental results for \cref{ex:4}} \label{fig:ex-complex-fun-grad}
\end{figure}

\begin{example}\label{ex:4}
We consider the JATD problem in \eqref{cost-fn-general-1} with $d=4$, $n_1=n_2=n_3=n_4=3$ and $r=2$. 
Let $L=1$ and $\alpha_1=1$.
Let $\tens{A}=\tens{A}^{\R}+\ui*\tens{A}^{\I}$, where $\tens{A}^{\R}$ and $\tens{A}^{\I}$ are generated randomly using $\textsf{round}(1 + \textsf{rand}(n_1,n_2,n_3,n_4))$, in which $\textsf{round}(\cdot)$ is the same as that in MATLAB.
It can be calculated that the initial ratio $\texttt{Rat}=0.0328$. In Jacobi-MG, we choose the parameter $\delta=\frac{1}{10}\sqrt{\frac{2}{d n_{\rm max}(n_{\rm max}-1)}}$.
The results of  $\|\tens{A}\|^2-f( \vect{\omega}_k)$ and $\|\ProjGrad{f}{\vect{\omega}_k}\|$ with respect to the running time are shown in \Cref{fig:ex-complex-fun-grad}. 
In particular, all the algorithms in \Cref{fig:ex-complex-fun-grad} arrive at the same ratio $\texttt{Rat}=0.9196$.
\end{example}

The results of the above examples demonstrate the efficiency and effectiveness of the tested Jacobi-type algorithms. 
Compared with the first-order Riemannian optimization methods,
the tested Jacobi-type algorithms have much better
performances in the convergence to stationary points, since the norm of the gradient for Jacobi-type algorithms is smaller from \Cref{fig:ex-joint-diag-fun-grad,fig:ex-complex-diag-fun-grad,fig:ex-complex-fun-grad}. On the other hand, in Jacobi-G, Jacobi-MG and Jacobi-GP algorithms, they need to compute the Riemannian gradient and check the gradient inequality \eqref{inequality-gra-based-2} or \eqref{inequality-gra-based} in each iteration, and thus cost more elapsed time. This is also consistent with our experimental results in  \Cref{fig:ex-joint-diag-fun-grad,fig:ex-complex-diag-fun-grad,fig:ex-complex-fun-grad}.

\section{Conclusion}\label{sec:conclusions}

 We considered two general classes of homogeneous polynomials with orthogonality constraints, the HPOSM problem \eqref{definition-f-new} and HPOSM-P problem \eqref{definition-f-prob-1}, which include as special cases the JATD-S and JATD problems. 
We proved some new convergence results for the Jacobi-G algorithm on the unitary group, which is to solve the JATD-S problem.  
We then proposed the Jacobi-MG algorithm on the product of unitary groups to solve the JATD problem, and also obtained several new convergence results.
Moreover, as a generalization of the Jacobi-PC algorithm proposed in \cite{LUC2017globally} and the JLROA-GP algorithm proposed in \cite{li2019jacobi}, we formulated the Jacobi-GP and Jacobi-MGP algorithms in the general case to solve the HPOSM problem \eqref{definition-f-new} and HPOSM-P problem \eqref{definition-f-prob-1}, and these two new proximal variants are proved to have global and weak convergence properties. 
In future work, it may be interesting to consider the following two questions: (i) how to establish the global convergence of Jacobi-G algorithm for JATD-S problem with complex matrices or order-$3$ complex tensors, or order-$4$ real symmetric tensors without any condition; (ii) how to solve the subproblem in Jacobi-G algorithm, if the subproblem can only be represented as the sum of a finite number of quadratic forms (not a single quadratic form as we assume in the present paper\fin).

\vspace{0.5cm}
\noindent 
{\bf Acknowledgments.} The authors would like to thank Prof. Shuzhong Zhang for his insightful discussions and suggestions.
 We are grateful to the associate editor and the two anonymous reviewers for their useful suggestions and comments, which helped to improve the presentation of the paper. 





\bibliographystyle{amsplain}
\bibliography{TensorRef}

\end{document}